\documentclass[noams]{compositio}

\renewcommand{\contentsname}{Contents\\{\footnotesize\normalfont(A table
of contents should normally not be included)}}

\usepackage{amsmath}
\usepackage{amssymb}
\usepackage{color}
\usepackage[colorlinks,anchorcolor=blue,citecolor=blue,linkcolor=blue,urlcolor =blue,bookmarksdepth=2]{hyperref}
\usepackage{multirow}
\usepackage{makecell}
\usepackage{mathtools}
\usepackage{comment}
\usepackage{tikz-cd}	
\usepackage{mathrsfs}   
\usepackage[flushleft]{threeparttable} 

\newtheorem{Thm}{Theorem}[section]
\newtheorem{Lem}[Thm]{Lemma}
\newtheorem{Cor}{Corollary}[Thm]
\newtheorem{Prop}[Thm]{Proposition}

\newtheorem{Def}{Definition}[section]
\newtheorem{eg}[Thm]{Example}

\newtheorem{Rmk}{Remark}[section]

\newenvironment{pf}{\begin{proof}}{\end{proof}}

\numberwithin{equation}{section}

\newcommand{\Z}{\mathbb{Z}}
\newcommand{\Q}{\mathbb{Q}}

\newcommand{\F}{\mathbb{F}}

\newcommand{\OO}{\mathcal{O}}

\newcommand{\RNum}[1]{\uppercase\expandafter{\romannumeral #1\relax}}

\newcommand{\Gal}{{\rm Gal}}

\newcommand{\GL}{{\rm GL}}

\DeclareFontFamily{U}{wncy}{}
\DeclareFontShape{U}{wncy}{m}{n}{<->wncyr10}{}
\DeclareSymbolFont{mcy}{U}{wncy}{m}{n}
\DeclareMathSymbol{\Sha}{\mathord}{mcy}{"58} 
\begin{document}

\title{Faltings-Serre method on three dimensional selfdual representations}
\author{Lian Duan}
\email{l.duanzwz@gmail.com}
\address{Department of Mathematics and Statistics, University of Massachusetts Amherst.
	710 N. Pleasant Street
	Amherst, MA 01003-9305, USA}
%
%

\begin{abstract}
We prove that a selfdual $GL_3$-Galois representation constructed by van Geemen and Top is isomorphic to a quadratic twist of the symmetric square of the Tate module of an elliptic curve. This is an application of our refinement of the Faltings-Serre method to $3$-dimensional Galois representations with the ground field not equal to $\Q$. The proof makes use of the Faltings-Serre method, $\ell$-adic Lie algebra, and Burnside groups.
\end{abstract}

\maketitle


\section{Introduction}

Let $K$ be a number field with ring of algebraic integers $\OO_K$ and absolute Galois group $G_K$. For each pair $(a,s)\in K^2$ with $a\neq \pm 1$ and $ s\neq 0$, consider the following elliptic surface parameterized by $t$
$$
y^2=x(x^2+2(s^2(a+1)+at^2)x+t^4)
$$
and take $\mathcal{S}_{a,s}$ to be its projective closure. van Geemen and Top \cite[$\S$~2.4]{GT-selfdual} consider a degree $4$ cover of $\mathcal{S}_{ a,s}$, denoted by ${\mathcal{A}}_{ a,s}$. They then construct a selfdual $3$-dimensional $G_K$ representation $(V_{\ell})_{ a,s}$ by taking a subquotient of the transcendental part of the second \'{e}tale cohomology of ${\mathcal{A}}_{a,s}$ \cite[$\S$5.1]{GT-selfdual}. Based on experimental data, they conjecture that for every such pair $(a,s)$, the representation $(V_{\ell})_{ a,s}$ is related to the symmetric square of the Tate module of an elliptic curve ${E}$ defined over $K$ or one of its quadratic extensions \cite[$\S 5.4$]{GT-selfdual}. Specifically, they conjecture that
\begin{equation}
\delta(N_{K/\Q}(\mathfrak{p}))  tr(F_{\mathfrak{p}}|(V_{\ell})_{ a,s})=tr(F_{\mathfrak{p}}|Sym^2(T_{{E}})),		\label{star}
\end{equation}
where $Sym^2(T_{{E}})$ is the symmetric square of the Tate module of ${E}$, $N_{K/\Q}$ is the norm, $tr$ is the trace, $F_{\mathfrak{p}}$ is the Frobenius class corresponding to the prime ideal $\mathfrak{p}$ of the algebraic integer ring of $K$, and $\delta$ is a Dirichlet character. 

Now let $K=\Q(\sqrt{-3})$, and let ${E}$ be the elliptic curve defined by 
$$
Y^2=X^3+(\sqrt{-1}-1)X^2+\left(-\frac{\sqrt{-1}}{4}+\frac{\sqrt{-3}}{8}-\frac{1}{8}\right)X.
$$
Take $V_{\ell}=(V_{\ell})_{\sqrt{-3}, 1}$, and write $\theta_{a}(*)=\left(\frac{a}{*}\right)$ for the Kronecker symbol corresponding to the integer $a$. In this paper we will prove the following theorem.
\begin{Thm}\label{Thm1}
	We have
	$$
	\theta_{-2}\otimes V_{\ell}\simeq Sym^2(T_{{E}}).
	$$
	In particular, equation $(\ref{star})$ is true for all $\mathfrak{p}$ not dividing $2$.
	
\end{Thm}

\begin{Rmk}
	Note that the elliptic curve ${E}$ is not defined over $K$. However it is a $K$-curve (cf. section \ref{Descend}), and thus its Tate module can be descended to a $G_K$-representation.
\end{Rmk}

Let $K$ be a number field, and let $M_{\lambda}$ be a finite extension of the $\ell$-adic field $\Q_{\ell}$. Let $\rho_1, \rho_2: G_K\to \GL_n(M_{\lambda})$ be Galois representations unramified outside a finite set of primes. In the proof of the Mordell Conjecture, Faltings shows that we can test if $\rho_1$, $\rho_2$ are equivalent up to semisimplification by performing a finite calculation \cite{Falting-FS}. Serre  \cite{Serre-FS} turns this into an effective tool, and Livn\'{e} \cite[thm.~4.3]{Livne-FS} improves this specifically for the case $n=2$. Many researchers (including but not limited to Boston \cite{Boston-Faltings-Serre}, Hulek, Kloosterman, Sch\"utt \cite{Schutt-Modularity-Calabi-Yau}, Schoen \cite{Schoen-85}, Socrates, Whitehouse \cite{Whitehouse-hilber-modular-form}, Dieulefait, Guerberoff, Pacetti \cite{Pacetti-Faltings-Serre-method}) have successfully applied the Faltings-Serre method to study the modularity over $\Q$ or some imaginary quadratic fields. However, due to the limits of current hardware, there was no known application of the Faltings-Serre method when $n>2$ and $K\neq \Q$.	

In his work \cite{Loic07}, Greni\'e finds an explicit bound of the norm of prime ideals such that the equivalence between $\rho_1$ and $\rho_2$ can be verified as long as they have the same traces for all Frobenius of unramified prime ideals with norm under this bound. However, a direct application of Greni\'e's work to the setups of Theorem  \ref{Thm1} leads to a bound that is too large to be verified. 

In this work, we refine Greni\'{e}'s criterion for $3$-dimensional selfdual Galois representations in two aspects. First, for general selfdual representations, by studying the rank of the Lie algebras of their images we reduce the number of prime ideals that needed to be checked in Greni\'e's result (see Theorem \ref{Thm2}). Second, suppose $K$ is quadratic, and suppose $\Gal(K^{ur}_{2,\infty}(2)/K)$ is generated by two elements, where $K^{ur}_{2,\infty}(2)$ is the maximal pro-$2$ extension of $K$ unramified outside $2$ and $\infty$. We improve the bound (Theorem \ref{Thm4}, Theorem \ref{Thm3-trick}) further by studying the structure of the Burnside group $B(2,4)$ (Example \ref{eg-burnside}). As an application, we verify Theorem \ref{Thm1} with the improved bound from Theorem \ref{Thm4}. In the following Table \ref{Table-Comparison}, we compare  Greni\'e's criterion, our first improvement (Theorem \ref{Thm2}) and our second improvement (Theorem \ref{Thm4}) in the case of Theorem \ref{Thm1}. In this table, we list the sizes of the sets $T$ of prime ideals that are needed to check to prove Theorem \ref{Thm1}, if the Extended Riemann Hypothesis (ERH) is assumed to find $T$, and the total time we spend to verify Theorem \ref{Thm1}.

\begin{table}[ht]
	\begin{center}
		\begin{threeparttable}
			\begin{tabular}{| c | c | c| c |}
				\hline
				Methods	&  size of $T$ & assume ERH & running time 	\\ \hline
				
				\makecell{By Greni\'e's criterion only\tnote{a}} 	& \multirow{2}{*}{unknown\tnote{a}} &  & \multirow{2}{*}{unknown} \\ \cline{1-1}
				
				\makecell{By first improvement only\tnote{a}} 	&  &&  \\ \hline
				
				\makecell{By Greni\'e's criterion and\\
					structure of $\Gal(K^{ur}_{2,\infty}(2)/K)$\tnote{b}} 	& \multirow{2}{*}{$\#(T)\geq 7\times 10^9$} & \multirow{2}{*}{yes} & \multirow{2}{*}{$>$ one year\tnote{c}} \\ \cline{1-1}
				
				\makecell{By first improvement and\\
					structure of $\Gal(K^{ur}_{2,\infty}(2)/K)$\tnote{b}} 	&  &  &  \\ \hline
				
				\makecell{Second improvement\\
					(Theorem \ref{Thm4})} 	& $\#(T)\leq 75$ & no & $\approx $ two weeks\tnote{d}	\\ \hline

			\end{tabular}
			
			\caption{Comparison of different methods}
			\label{Table-Comparison}
			\begin{tablenotes}
				\small 
				\item[a] In these cases, then the field $K_3$ (see below) is too large to be constructed by computers. 
				\item[b] Details can be found in section \ref{Find-covering-set}.
				\item[c] This is an estimation, in fact we did not finish $10 \%$ of the process after two months. 
				\item[d] We spent about two weeks finding $T$. Once $T$ is found, it takes less than one day to verify Theorem \ref{Thm1}.
			\end{tablenotes}
		\end{threeparttable}		
	\end{center}
	
\end{table}

To state the first refinement, we introduce the following concept and notations. A matrix in $\GL_n(\Z_{\ell})$ is \textit{congruent trivial} if its characteristic polynomial is congruent to $(t-1)^n$  $(\bmod \ell)$. We say a $G_K$-representation is \textit{congruent trivial} if every element in the image is congruent trivial. Take $S$ to be a finite set of prime places ($S$ may include the Archmedian places). For a pair of $\ell$-adic $G_K$-representations $(\rho_1,\rho_2)$, let $K_0=K$, then let $K_1$ to be the Galois extension of $K_0$ such that $K_1$ is unramified outside $S$ and $\Gal(K_1/K_0)$ is isomorphic to the image of the $\bmod\ \ell$ residue representation of $\rho_1\oplus \rho_2$. Then for each $i\geq 1$, let $K_{i+1}$ be the maximal abelian extension over $K_i$ unramified outside  $S$ and $\Gal(K_{i+1}/K_i)$ is an elementary $\ell$-group. Let $\epsilon=1$ if $\ell=2$ or $0$ otherwise, then let $K_S=K_{5+3\epsilon}$. Take $T$ to be a finite set of prime ideals in $\OO_K$ such that $\mathfrak{p}\nmid 2$ and every element of $\Gal(K_S/K)$ corresponds to at least one Frobenius $F_{\mathfrak{p}}$ with $\mathfrak{p}\in T$ (i.e. $T$ is a covering set of $\Gal(K_S/K)$, see Definition \ref{covering-set}). For an $\ell$-adic Galois representation $\varphi$, its $m$th Tate twist $\varphi\otimes \mu_{\ell}^{\otimes m}$ is denoted by $\varphi(m)$, where $\mu_{\ell}$ is the cyclotomic representation. Denote by $\varphi^*$ the dual representation of $\varphi$.

\begin{Thm}\label{Thm2}
	With notations as above, suppose that $\rho_1, \rho_2:G_K\to GL_3(\Z_{\ell})$ both satisfy $\rho_i^*\simeq \rho_i (2m)$ for some integer $m$, and suppose that both $\rho_i(-m) $ are congruent trivial. Moreover, suppose $T$ is disjoint with the ramified ideals with respect to $\rho_1\oplus \rho_2$. Then $\rho_1$ and $\rho_2$ are equivalent up to semisimplification if 
	$$
	tr(\rho_1(F_{\mathfrak{p}}))=tr(\rho_2(F_{\mathfrak{p}}))
	$$
	for all $\mathfrak{p}\in T$.
\end{Thm}

\begin{Rmk}\label{Compracsion1}
	As a comparison, using Greni\'e's criterion, we have to take $K_S$ to be $K_{8+3\epsilon}$. With our criterion, $K_S=K_{5+3\epsilon}$. For the purpose of Theorem \ref{Thm1}, our criterion reduces the degree of $K_S$ by a factor of at least $2^9$.
\end{Rmk}

Theorem \ref{Thm2} is not enough to prove Theorem \ref{Thm1} since it takes more than a month to construct a degree $2^7$ extension using PARI/Magma while $K_2$ has degree $2^8$. In order to prove Theorem \ref{Thm1}, we reduce the size of the set $T$ further in Theorem \ref{Thm2} by making use of the fact that $\Gal(K^{ur}_{2,\infty}(2)/K)$ is a free pro-$2$ group with two generators \cite[Theorem~2]{Jossey-2-ext-ur-outside-2} and studying the Burnside group $B(2,4)$. 

\begin{Thm}\label{Thm4}
	Let $K=\Q(\sqrt{n})$ with $n=-1,-2,-p$ or $-2p$ where $p=\pm 3\ (mod\ 8)$. Assume $\rho_i$ $(i=1,2)$ are congruent trivial and unramified outside $2$ and $\infty$. Then there exists a set $T$ which only depends on $K$ and consists of at most $75$ prime ideals of $\OO_K$ not lying above $2$. With this $T$, we have that $\rho_1$ and $\rho_2$ are equivalent if and only if  
	$$
	tr(\rho_1(F_{\mathfrak{p}}))=tr(\rho_2(F_{\mathfrak{p}}))
	$$
	for all $\mathfrak{p}\in T$. In particular, when $K=\Q(\sqrt{-2})$, the set $T$ is given by Table \ref{List2}. When $K=\Q(\sqrt{-3})$, $T$ is given by Table \ref{List3}.
\end{Thm}

\begin{Rmk}
	Theorem \ref{Thm4} is effective in the sense that all elements in the finite set $T$ can be listed by Theorem \ref{Thm3-trick}.
\end{Rmk}

\begin{Rmk}\label{Rmk-nonselfdual-comparison}
	Theorem \ref{Thm4} and hence Theorem \ref{Thm3-trick} also work for non-selfdual representations except that when comparing non-selfdual representations, one needs to check 
	$$
	char(\rho_1(F_{\mathfrak{p}}))=char(\rho_2(F_{\mathfrak{p}}))
	$$
	for all $\mathfrak{p}\in T$. Here $char$ stands for the characteristic polynomial. 
\end{Rmk}

One can see that Theorem \ref{Thm1} immediately follows once we verify that the two representations $V_{\ell}$ and $Sym^2(T_E)$ are both congruent trivial.

Here is an outline of this paper. In section \ref{section-background}, \ref{Section-FS-method} and \ref{section-pro-p} we review the background of Galois representations, Faltings-Serre method and pro-$p$ groups respectively. In particular, we will review powerful pro-$p$ groups, and recall a theorem that will help us to find a powerful subgroup in every pro-$p$ group.

We prove our first improvement in section \ref{section-Lie-alg}. Given a $3$-dimensional selfdual representation $\rho$, in order to find a bound of the rank of its image, we study the Lie algebras of its image. We show the dimension of its Lie algebra is at most $3$, and therefore the rank of its image is at most $3$. Hence Theorem \ref{Thm2} follows as a consequence.

In section \ref{section-pf-Thm1}, we prove Theorem \ref{Thm1}. First, we descend the Tate module of the elliptic curve $E$ in Theorem \ref{Thm1} as a $G_K$-representation by a result of Ribet. Also, we give the formula to compute the trace of the symmetric square of the descended representation. Then to compare the two sides of $($\ref{star}$)$, we try to find a covering set $T$ by Theorem \ref{Thm2}. Then to speed up this process, we prove Theorem \ref{Thm3-trick} and successfully cut off the size of $T$ so that the whole process can be completed in two weeks. Hence Theorem \ref{Thm1} is verified.

\begin{acknowledgements}
We would like to thank Professor Siman Wong for suggesting this exciting topic and giving helpful advice. We thank Professor Paul Gunnells and Professor Farshid Hajir for conversation and suggestions in representation theory and background of pro-$p$ groups. We thank Professor Lo\"ic Greni\'e for helpful suggestion in computation and thank Professor Hans Johnston for providing computing source. 
\end{acknowledgements}

\subsection*{Notations} In this paper, unless mentioned specifically, we will assume the followings.
\begin{itemize}
	\item $\Q$ is the rational field, with integer ring $\Z$ and $\ell, p$ represent prime integers of $\Z$.
	\item $K, L$ represent number fields or infinitely Galois extensions of $\Q$, with algebraic integer ring $\OO_K, \OO_L$ respectively. Prime ideals are usually written as $\mathfrak{p}$ or $\mathfrak{P}$. Their corresponding Frobenius class are denoted by $F_{\mathfrak{p}}$ or $\mathfrak{P}$ respectively. 
	\item Given a number field $K$ and a finite set $S$ of prime places of $K$, $K^{ur}_{S}$ is the maximal Galois extension of $K$ which is unramified outside $S$. 
	\item $\Q_{\ell}$ is the local field with ring of integer $\Z_{\ell}$. $M_{\lambda}$ is either a finite extension of $\Q_{\ell}$ or the algebraic closure $\overline{\Q}_{\ell}$, with $\mathcal{O}_{M_{\lambda}}$ (or $\mathcal{O}$ is there is no confusion) its integer ring. 
	\item If $F'/F$ is an Galois extension of either local or global field, we denote by $\Gal(F'/F)$ the corresponding Galois group. In particular, $G_{F}=\Gal(\overline{F}/F)$ is the absolute Galois group of $F$.
	\item $K^{ur}_{S}(p)$ is the maximal pro-$p$ extension of $K$ unramified outside $S$. In particular, $K^{ur}_{2,\infty}(2)$ is the maximal pro-$2$ extension of $K$ unramified outside $2$ and $\infty$.
	\item $\rho$, $\varphi$ represent Galois representations. 
	\item $\mathcal{G}$ represents a Lie group with its Lie algebra $\mathfrak{g}$.
	\item $E$, $\tilde{E}$ are elliptic curves defined over number fields. Their Tate module are written as $T_{\ell}(E)$ and $T_{\ell}{\tilde{E}}$ respectively. Without danger of confusion we will simply write $T_{\ell}$.
	\item $T$ stands for a covering set (Definition \ref{covering-set}).
\end{itemize}

\section{Background of Galois representations}\label{section-background}

In this section, we recall the background of $\ell$-adic Galois representations. We denote by $\Gal(L/K)$ the Galois group for the (finite or infinite) Galois extension $L/K$, equipped with the pro-finite topology. In particular, we write $G_K$ for $\Gal(\overline{K}/K)$. With these notations, an (n-dimensional) $\ell$-adic representation of $\Gal(L/K)$ is a continuous group homomorphism 
$$
\rho: \Gal(L/K)\to GL_n(M_{\lambda}).
$$
Since $\Gal(L/K)$ is a compact group,  we can assume that their images are in $GL_n(\mathcal{O})$ (cf. \cite[Prop. 9.3.5]{Diamond-GTM-228} or \cite[Remark~1, p.~\RNum{1}-1]{Serre-l-adic-repn}), where $\OO$ is the ring of integers of $M_{\lambda}$.

Fix a number field $K$, and let $\mathfrak{p}$ be one of the prime ideals of its algebraic integer ring $\mathcal{O}_K$. Then we have corresponding local field $K_{\mathfrak{p}}$ and the corresponding residue field $k_{\mathfrak{p}}$. The kernel of the natural quotient $G_{K_{\mathfrak{p}}}\to G_{k_{\mathfrak{p}}}$ is the inertia group at $\mathfrak{p}$, and denoted by $I_{\mathfrak{p}}$, and we denote by $F_{\mathfrak{p}}$ the preimage of the Frobenius of $G_{k_{\mathfrak{p}}}$. Take $\rho$ to be a $\Gal(L/K)$-representation. It is called \textit{unramified} at $\mathfrak{p}$ if $I_{\mathfrak{p}}$ is in the kernel of $\rho$. Moreover, when $\rho$ is unramified at $\mathfrak{p}$, it makes sense to consider $\rho(F_\mathfrak{p})$ as an element in $GL_n(M_{\lambda})$. Let $S$ be a set consisting of finite prime ideals (which may include the infinite primes). We denote by $K^{ur}_S$ the maximal Galois extension above $K$ which is unramified outside $S$. Thus if $\rho$ is unramified outside $S$ then it factors through $\Gal(K^{ur}_S/K)$.

Recall that every representation has a Jordan-H\"older composition series. Two representations $\rho_1$, $\rho_2$ are said to be \textit{equivalent} (up to semisimple) if they share the same Jordan-H\"older composition series up to a reorder. In this case, we write $\rho_1\sim \rho_2$. A representation $\rho$ is called \textit{semisimple} if it is a direct sum of simple representations, or equivalently, if every proper sub-representation of $\rho$ has its complement as a sub-representation. According to the Jacobian density theorem \cite[chap. 4, $\S$ 3]{Jacobson-BA2}, the equivalent class of an $\ell$-adic Galois representation is uniquely determined by its traces.

\begin{Prop}\label{Prop-2}
	Let $\rho_i : G\to GL_n(M_{\lambda})$ ($i=1,2$) be two $\ell$-adic Galois representations. Then
	$
	\rho_1\sim \rho_2
	$
	if and only if 
	$$
	tr(\rho_1(g))=tr(\rho_2(g))
	$$
	for all $g\in G$. Here $tr$ means the trace of a representation. 
\end{Prop}

Given a representation $\rho$, there is a natural \emph{dual representation} $\rho^*$. Precisely, it can be computed by $\rho^*(g)=(\rho(g)^{T})^{-1}$, i.e. the transverse inverse of $\rho(g)$.

\begin{Def}\label{sefldual-Def}
	A $G_K$-representation $\rho$ is called \textit{selfdual} if $\rho^*$ is isomorphic to $\rho(2m)$ for some integer number $m$. Here $\rho(2m):=\rho \otimes \mu_{\ell}^{\otimes 2m}$ is a Tate twist of $\rho$, where $\mu_{\ell}$ is the cyclotomic representation. In particular, if $m=0$ i.e. $\rho\simeq \rho^*$, then we say that $\rho$ is \textit{strictly selfdual}.
\end{Def}

In the rest of this section, we focus on two kinds of selfdual representations related to Theorem \ref{Thm1}. 
\subsection{Symmetric square of Tate module of elliptic curves}

Given an elliptic curve ${E}$ over $K$, and let $\ell$ be a prime integer. The corresponding Tate module $T_{\ell}(E)$ induces a $2$-dimensional $G_K$ representation, denoted by $\varphi$. If $\mathfrak{p}$ does not divide the discriminant of ${E}$, then $\varphi$ is unramified at $\mathfrak{p}$ and $tr(\varphi(F_{\mathfrak{p}}))=\alpha+\beta\in \Q$ for some $\alpha, \beta\in \overline{\Q}$ which are independent of $\ell$. Both of $\alpha$ and $\beta$ have absolute value $\sqrt{q}$ with $q={N_{K/\Q}(\mathfrak{p})}$, and $\alpha \beta=q$. And one can see that $tr(F^{-1}_{\mathfrak{p}})={\alpha}/{q}+{\beta}/{q}$.

Now let $\rho$ be the symmetric square of $\varphi$. A simple calculation shows that 
$$
tr(\rho(F_{\mathfrak{p}}))=tr(\varphi(F_{\mathfrak{p}}))^2-q=(\alpha+\beta)^2-q
$$
and 
$$
tr(\rho(F^{-1}_{\mathfrak{p}}))=({\alpha}/{q}+{\beta}/{q})^2-{1}/{q}.
$$
Thus, $tr(\rho^*(2)(F_{\mathfrak{p}}))=tr(\rho^*(F_{\mathfrak{p}}))q^2=(({\alpha}/{q}+{\beta}/{q})^2-{1}/{q})q^2=tr(\rho(F_{\mathfrak{p}}))$ for all unramified $\mathfrak{p}$. Then by the Chebotarev density theorem we know that $\rho$ is selfdual up to semisimplification. 

For later use, we cite the Serre's open image theorem.
\begin{Prop}\cite[$\S$2.2 Theorem at page \RNum{4}-12]{Serre-l-adic-repn}\label{Thm-Serre-open-image}
	Let ${E}$ be an elliptic curve defined over a number field $K$, and let $\varphi: G_K\to GL_2(\Q_{\ell})$ be the Galois representation induced by the Tate module of ${E}$. If ${E}$ has no complex multiplication, then $img(\varphi)$ has finite index in $GL_2(\Z_{\ell})$.
\end{Prop}

\begin{Cor}\label{Cor-Sym-irr}
	Let $V$ be the representation space of $Sym^2(\varphi)$. Then $V\otimes \Q_{\ell}$ is generated by any nonzero vector as a $\Q_{\ell}[G_K]$-module. In particular, $V$ is an irreducible $G_K$-representation. 
\end{Cor}
\begin{pf}
	This follows from Proposition \ref{Thm-Serre-open-image} and direct calculation. 
\end{pf}

\subsection{The selfdual representation of van Geemen and Top}

In their paper \cite[$\S$ 2]{GT-selfdual}, for each $\mathcal{S}_{ a,s}$, van Geemen and Top construct a degree $4$ branched covering surface $\mathcal{A}_{ a,s}$, i.e., there is a degree $4$ automorphism $\sigma$ of $\mathcal{A}_{ a,s}$ such that $\mathcal{A}_{ a,s}/\langle \sigma \rangle=$ $\mathcal{S}_{ a,s}$. They consider a subquotient of the transcendental part of the second \'etale cohomology of $\mathcal{A}_{ a,s}$, and find that it is a Galois representation which admits an action induced by $\sigma$. Then the representation space $(V_{\ell})_{ a,s}$ is defined to be one of the eigenspaces of $\sigma$. They show that when $a\neq \pm 1$ and $ s\neq 0$, the corresponding $G_K$-representation on $(V_{\ell})_{ a,s}$ is $3$-dimensional and (possibly up to semisimplification) selfdual \cite[Prop.~5.2]{GT-selfdual}. In that case, they show that \cite[Prop.~3.1 and Thm.~3.5]{GT-selfdual}
$$
tr(F_{\mathfrak{p}}|(V_{\ell})_{ a,s})=\#(\mathcal{S}_{ a,s})_{\infty}(\F_q)+\sum\limits_{\tau\in \F_q, \tau^2+16\in \F_q^{*2}} \left( \frac{u_\tau^2+4}{q} \right) \#(\mathcal{S}_{ a,s})_{r+s\tau/4}(\F_q)
$$
where $\F_q$ is the residue field corresponding to finite prime ideal $\mathfrak{p}$ and $u_\tau$ denotes a root in $\F_q$ of $X^2-\tau X-4=0$ and $(\mathcal{S}_{ a,s})_{t}$ is the fiber over $t$.

Specifically, when $( a,s)=(\sqrt{-3},0,1)$, the surface is $$
\mathcal{S}:Y^2=X(X^2+2(\sqrt{-3}+1+\sqrt{-3}t^2)X+t^4).
$$

When $\ell=2$, $V_{\ell}$ is defined over $\Z_2(\sqrt{-1})$. But its semisimplification is in fact a $GL_3(\Z_2)$ representation according to the following proposition and the fact that its characteristic polynomial of $F_{\mathfrak{p}}$ with prime $\mathfrak{p}$ above $31$ has three distinct roots in $\Z_2$.  

\begin{Prop}
	Let $G$ be a group, and $E$ a field of characteristic $0$. Let $\phi: G\to GL_n(\overline{E})$ be a semisimple representation defined over $\overline{E}$. let $\phi\simeq \phi_1\oplus \cdots \phi_r$ be an irreducible decomposition of $\phi$. Assume that the following conditions are satisfied:
	\begin{enumerate}
		\item $\phi$ is defined over a finite extension of $E$.
		\item $tr(\phi)\in E$ for every $g\in G$.
		\item There is an elemetn $g_0\in G$ such that the characteristic polynomial of $\phi(g_0)$ has $n$ distinct roots in $E$.
	\end{enumerate}
	Then each $\phi_i$ is defined over $E$. In particular, $\phi$ is defined over $E$.
\end{Prop}

\begin{pf}
	See \cite[Prop. 7]{Chin-Galois-repn-descend}, \cite[Lem. 2.1]{Ito-K-M-non-selfdual} or the proof of \cite[Prop. 3.2.5]{Chenevier-Harris-construction-auto-repn}.
\end{pf}

Fix $\mathfrak{p}\nmid 2$, $V_{\ell}$ is unramified at $\mathfrak{p}$, and we claim that the characteristic polynomial of $F_{\mathfrak{p}}$ with respect to $V_{\ell}$ satisfies 
\begin{equation} \label{eqn-1}
t^3+t^2+t+1=(t-1)^3\ \ (mod\ 2).
\end{equation}

To see this, note that $V_{\ell}$ is a selfdual representation, so we only need to compute its trace of $F_{\mathfrak{p}}$. Making use of the above formula for trace, and the symmetry of $S$ with respect to the involution $t\mapsto -t$, and also reviewing the details of the construction of $V_{\ell}$ \cite[$\S$.2]{GT-selfdual}, we can compute $\#\mathcal{S}_{t}(\F_q)$ $(mod\ 4)$ as $t$ runs through $\mathbb{P}_{\F_{q}}^1$. In fact, we have the following:
\begin{enumerate}
	\item[(a)] When $t=0$, if $\sqrt{3}+\sqrt{-1}\in \F_q$, then $\mathcal{S}_0$ contributes $q$ points; otherwise it contributes $q+2$ points.
	
	\item[(b)] When $t=\infty$, we have $\# \mathcal{S}_{\infty}(\F_q)=0\ (mod\ 4)$. 
	
	\item[(c)] When $t=\pm \frac{1+\sqrt{-3}}{2}$, if $\sqrt{-1}\in \F_q$, then $\mathcal{S}_t$ contributes $q$ points; otherwise it contributes $q+2$ points. 
	
	\item[(d)] When $t=\pm i$, if $\sqrt{-1}\in \F_q$, then $\mathcal{S}_t$ gives $q$ points; otherwise it contributes $q+2$ points. 
\end{enumerate}  
By all above, we obtain
$$
tr(F_{\mathfrak{p}}|V_{\ell})=
\begin{cases}
1 (\ mod\ 4)\text{ if }p\equiv 7,13\ (mod\ 24)\\
3 (\ mod\ 4)\text{ if }p\equiv 1,5,11,17,19,23\ (mod\ 24)
\end{cases}
$$
and our claim follows immediately.

\section{Faltings-Serre method}\label{Section-FS-method}
In this section, we recall the Faltings-Serre method. Some useful references are \cite[chap.~5, $\S$.~5.2, 5.4]{Gabriell-thesis} and \cite{Livne-FS}. We will follow \cite{Gabriell-thesis}. If $S$ to be a finite set of prime places of number field $K$, and we assume that $\rho_i: G_K\to \GL_n(\mathcal{O}_{\lambda})$ $(i=1,2)$ are unramified outside $S$. Take $\rho=\rho_1\oplus \rho_2$, and consider it as an $\mathcal{O}_{\lambda}$-algebra homomorphism
$$
\rho: \mathcal{O}_{\lambda}[G_K]\to M_n( \mathcal{O}_{\lambda})\oplus M_n(\mathcal{O}_{\lambda}).
$$
Let $N$ be its image and consider the composition 
$$
\delta: G_K\to N^{\times}\to (N/\lambda N)^{\times}.
$$

\begin{Def}
	The image $\delta(G_K)$ is called the \textit{deviation group} of the pair $(\rho_1, \rho_2)$.
\end{Def}

\begin{Rmk}
	$\delta(G_K)$ a finite group. But in general it is not a subgroup of $\GL_n(k)\oplus \GL_n(k)$ \cite[Proposition~5.2.2 and its remark]{Gabriell-thesis}.
\end{Rmk}

The following proposition improves Proposition \ref{Prop-2}.
\begin{Prop}\cite[Prop.~5.2.3]{Gabriell-thesis}\label{trace-Sigma}
	Let $\Sigma$ be a subset of $G_K$ surjecting onto $\delta(G_K)$. Then 
	$$
	\rho_1\sim \rho_2 \Leftrightarrow tr\rho_1|_{\Sigma}=tr \rho_2|_{\Sigma}.
	$$
\end{Prop}

\begin{Def}\label{covering-set}
	Fix a number field $K$, let $U$ be a finite set, and $\Psi: G_K\to U$ be a map of sets. A finite set $T$ of prime ideals of $\OO_K$ is called a \textit{covering set} of $U$ (with respect to $\Psi$) if every element of $U$ is in the image $\Psi(F_{\mathfrak{p}})$ for at least one $\mathfrak{p}\in T$.
	
	In particular, if $L$ is a finite Galois extension of $K$, $T$ is called a \textit{covering set of $\Gal(L/K)$} if it is a covering set of $\Psi: G_K\to U$ with $U=\Gal(L/K)$ and with $\Psi$ to be the natural quotient map $G_K\to \Gal(L/K)$.
\end{Def}

\begin{Rmk}\label{Rmk-covering-set}
	Using this definition, we can restate Proposition \ref{trace-Sigma} as follows:
	$$\rho_1\sim \rho_2 \textit{ if and only if } tr\rho_1(F_{\mathfrak{p}})=tr\rho_2(F_{\mathfrak{p}}) \textit{ for all }\mathfrak{p}\in T$$
	\textit{where }$T$\textit{ is a covering set of }$\delta(G_K)$. 
	
	\textit{In particular, if $\rho_1$ and $\rho_2$ can be descended to $\Gal(L/K)$-representations, with $L/K$ a finite Galois extension, then  $T$ can be chosen as any covering set of $\Gal(L/K)$}. 
\end{Rmk}

For the rest of this section, to simplify our arguments, we will always assume that $n=3$, and $M_{\lambda}=\Q_{\ell}$ and we assume the following congruent trivial condition for $\rho_i$ $(i=1,2)$:
\begin{equation}\label{cong-trivial-assumption}
\{\text{the characteristic polynomial of }\rho_i(g)\}\equiv (t-1)^3\ (\bmod \ell) \text{ for all }g\in G_K.
\end{equation}

\begin{Prop}
	Under the assumption $(\ref{cong-trivial-assumption})$, $\delta(G_K)$ is an $\ell$-group.
\end{Prop}
\begin{pf}
	In fact, we have a filtration of the image of $\rho$:
	$$
	img(\rho(G_K))=:G_0>G_1>G_2>\cdots>G_m>\cdots
	$$
	where $G_i$ for $i\geq 1$ is the kernel of $img(\rho(G))\ (mod\ \lambda^i)$. Since for every $i\geq 1$, the quotient $G_i/G_{i+1}$ is isomorphic to a subgroup of $(\Z/\ell\Z)^{18}$, hence $G_1$ is pro-$\ell$. Now consider $G_0/G_1$ as a subgroup of $GL_3(k)\oplus GL_3(k)$. When $\ell\geq 3$, let $t\in G_1/G_0$ we have 
	$$t^{\ell}-1=(t-1)^{\ell}=(t-1)^3\times (t-1)^{\ell-3}=0$$
	which means every element in $G_1/G_0$ has order $\ell$. When $\ell=2$, by listing all possible characteristic polynomials of elements in $GL_3(\F_2)$ we can show the same result. 
\end{pf}

\begin{Def}
	With $G_K$, $\rho_1, \rho_2$ as above, define $\Xi$ to be the subset of elements $g\in G_K$ for which the characteristic polynomials of $\rho_1(g)$ and $\rho_2(g)$ coincide (or, equivalently, $tr_1(g^i)=tr_2(g^i)$ for $i=1,2,3$).
\end{Def}
By Proposition \ref{trace-Sigma} (and Remark \ref{Rmk-covering-set}), we know that if $\Xi$ surjects onto $\delta(G_K)$, then $\rho_1\sim \rho_2$. Thus showing that $\rho_1\sim \rho_2$ is reduced to finding at least one subset of $\Xi$ which is a covering set of $\delta(G_K)$.

\begin{Prop}\label{exp4Prop}
	If $g\in \Xi$, then $\delta(g)^{\ell}=1$ $(\delta(g)^4=1$ when $\ell=2)$ in $\delta(G_K)$.
\end{Prop}

\begin{pf}
	Recall that $g\in \Xi$ means that $\rho_1(g)$ and $\rho_2(g)$ have the same characteristic polynomial. Denoteing this common polynomial by $t^3+a_2t^2+a_1t+a_0$. If $\ell>2$, we know that 
	$$
	(t-1)^{\ell-3}(t^3+a_2t^2+a_1t+a_0)=(t-1)^{\ell-3}(t-1)^3=(t-1)^{\ell}=t^{\ell}-1\ (\bmod \ell).
	$$ 
	Thus one can see that $\rho(g)^{\ell}-1\in \lambda(M)$, which implies that  $\delta(g)^{\ell}=1$ in $\delta(G_K)$. By similar argument we could show that we get the corresponding result for the case $\ell=2$.
\end{pf}

\begin{Def}\label{Def-G_l}
	Given a group $G$, denote by $G[\ell]$ (resp.~$G[4]$ if $\ell=2$) the subset of elements of order dividing $\ell$ (resp.~$4$) in $G$, and let $G^{\ell}=\langle g^{\ell}|h\in G\rangle$ (resp.~$G^{4}=\langle g^{4}|h\in G\rangle$) be the subgroup generated by the $\ell$th (resp.~$4$th) power of elements in $G$, and let $G_{\ell}=G/G^{\ell}$ (resp.~$G_4=G/G^4$).
\end{Def}

\begin{Lem}\label{ltorsion-lemma}
	Let $H$ be a $(pro)$ $\ell$-group such that every element in $H_{\ell}$ (resp.~ $H_4$ when $\ell=2$) has a lift to an element of $H[\ell]$ (resp.~$H[4]$). Then $H=H_{\ell}$ (resp.~$H=H_4$).
\end{Lem}

\begin{pf}
	\cite[Lemma.~7]{Loic07} or \cite[Lemma.~5.4.7]{Gabriell-thesis}).
\end{pf}

\begin{Prop}\label{Prop-compare-repns}
	If $\rho_1$ and $\rho_2$ be two representations satisfying the condition $(\ref{cong-trivial-assumption})$. Then the followings are equivalent:
	
	\begin{enumerate}
		\item $\rho_1\sim \rho_2$;
		\item $\Xi$ is a covering set of $(G_{K})_{\ell}$ $($resp.~$(G_{K})_4$$)$;
		\item $\Xi$ is a covering set of $\rho(G_K)_{\ell}$ $($resp.~$\rho(G_K)_4$$)$.
	\end{enumerate} 
\end{Prop}

\begin{pf}
	The only non-trivial part is $(3)\Rightarrow(1)$. To prove this implication, we first take $\ell>2$, and consider the following commutative diagram
\begin{center}
    \begin{tikzcd}
        G_K \arrow[r] \arrow[d]& \rho(G_K)_{\ell} \arrow[d]\\
        \delta(G_K) \arrow[r]& \delta(G_K)_{\ell}
    \end{tikzcd}
\end{center}

	Since $\Xi$ covers $(G_K)_{\ell}$, it also covers $\delta(G_K)_{\ell}$. Every element $\bar{g}$ in $\delta(G_K)_{\ell}$ has a lifting to $\Xi$, denoted by $g$, and by Proposition \ref{exp4Prop}, we know that $\delta(g)\in \delta(G_K)[\ell]$. Then the conclusion follows from the Lemma \ref{ltorsion-lemma}. By a similar argument we also show the result for case $\ell=2$.
\end{pf}

\section{Background of (pro)-$p$ groups}\label{section-pro-p}
In this section, we collect the necessary background on (pro)-$p$ groups, which will be helpful later. In particular, we will use Theorem \ref{powerful subgp} in the next section. We adopt the definitions and notations in \cite{Dixon-pro-p}, especially its first three chapters. 

Let $p$ be a prime integer. A group is called a $p$-group if each element of this group has a power of $p$ as its order. A pro-$p$ group is a topological profinite $p$-group, or equivalently, it is an inverse limit of a system of finite $p$-groups with respect to profinite topology. For the rest of this section, unless otherwise stated, all groups are assumed to be (pro)-$p$. For two elements $g$ and $h$ in group $G$, we denote $[g,h]$ to be  the commutator of $ghg^{-1}h^{-1}$, and  $[G,G]:=\{[g,h]| g,h\in G\}$, and $G^p:=\langle g^p|g\in G\rangle $ is the normal subgroup generated by the $p$th power of elements in $G$.

\begin{Def}
	For a (pro)-$p$ group $G$, the \textit{Frattini subgroup} of $G$, denoted by $\Phi(G)$, is the intersection of all maximal proper subgroups of $G$.
\end{Def}

\begin{Prop}\cite[0.9]{Dixon-pro-p} \label{Prop-3}
	\begin{enumerate}
		\item $\Phi(G)=[G,G]G^p$.
		\item Let $X\subset G$ be a subset, and assume that $X\Phi(G)$ generates $G/\Phi(G)$. Then $X$ generates $G$.  
		\item Let $d$ be the minimal cardinality of any topological generating set for $G$. Then $G/\Phi(G)\simeq \F_p^d$, and we denote by $d(G)$ the number $d$.
	\end{enumerate}
\end{Prop}

\begin{Def}
	For a finite $p$-group $G$, we define the rank of $G$ to be
	$$
	rank(G):=\sup\{d(H)|H<G\}.
	$$
\end{Def}

\begin{Def}\cite[3.11]{Dixon-pro-p}\label{Def-rank-pro-p}
	For a pro-$p$ group $G$, we define the rank of $G$ to be the common value of following $r_i$ ($i=1,2,3,4$):
	\begin{align*}
	r_1 &= \sup \{d(H)| H \text{ is a closed subgroup of }G \}\\
	r_2 &= \sup \{d(H)| H \text{ is a closed subgroup of } G \text{ and }d(H)<\infty \}\\
	r_3 &= \sup \{d(H)| H \text{ is a open subgroup of }G \}\\
	r_4 &= \sup \{rank(G/N)| N \text{ is a open subgroup of }G \}
	\end{align*}
\end{Def}

\begin{Def}
	The \textit{exponent} of a group  (not necessary profinite) $G$ is the least common multiple of order of elements in $G$. 
\end{Def}

\begin{eg}
	By definition of $G^p$, we know that $G/G^p$ has exponent $p$. In fact, it is the largest quotient of $G$ with this property, i.e., every other exponent $p$ quotient of $G$ has to factor through $G/G^p$.
\end{eg}

\begin{eg}\label{eg-burnside}
	Given a free (and not a profinite) group $F$ generated by $d$ elements and $n$ a positive integer, we denote by $B(d,n)$ the quotient $F/F^n$. It is called the \textit{Burnside group}. For a fixed pair $(d,n)$, the Burnside group is the universal group which is generated by $d$ elements and has exponent $n$. For all $d$, the group $B(d,2)$ is isomorphic to $(\Z/2\Z)^{\oplus d}$. However, when $n>2$, not much is known about $B(d,n)$. In this article, we will use the groups structure of $B(2,4)$ \cite{Tobin-Burnside-gp}.
\end{eg}

The definition of powerful (pro)-$p$ groups and its properities will help us to find a subgroup of $G^p$ ($G^4$ if $p=2$) for a given $G$.
\begin{Def}\label{Def-powerfun-subgp}
	A (pro)-$p$ group $G$ is called \textit{powerful} if $G/G^p$ (or $G/G^4$ when $p=2$) is abelian.
\end{Def}

\begin{Cor}\cite[3.4]{Dixon-pro-p}\label{powerfulcor}
	Let $G$ be powerful and finitely generated, then every element of $G^p$ is a $p$th power in $G$ and $G^p$ (resp. $G^4$) is open in $G$.  
\end{Cor}
\begin{pf}
	This is an immediate consequence of Definition \ref{Def-powerfun-subgp}. 
\end{pf}

The following theorem will help us in finding a powerful subgroup for a given $G$. 
\begin{Thm}\cite{Dixon-pro-p, Loic07}\label{powerful subgp}
	Let $G$ be a finitely generated pro-$p$ group of rank $r$, and define $\lambda(r)$ be the minimal integer such that $2^{\lambda(r)}\geq r$. Then $G$ has a powerful open subgroup of index at most $p^{r\lambda(r)}$ if $p$ is odd, and $2^{r+r\lambda(r)}$ if $p=2$. 
\end{Thm}

\begin{Rmk}\label{Rmk-filtration}
	In fact, the method of the proof is to construct a filtration of subgroups 
	$$
	G=:G_0>G_1>\cdots >G_{\lambda(r)+\epsilon}=:V
	$$
	such that $V$ is a powerful subgroup of $G$, where $\epsilon=1$ when $\ell=2$ and $0$ otherwise.
\end{Rmk}

\section{Selfdual Lie algebras and proof of Theorem \ref{Thm2}}\label{section-Lie-alg}

From the first four sections, especially by Proposition \ref{trace-Sigma}, we know that to compare $\rho_1$ and $\rho_2$ that satisfy the condition (\ref{cong-trivial-assumption}), we need to find a covering set of $G_{\ell}$ (resp. $G_4$ is $\ell=2$), where $G$ is $(\rho_1\oplus \rho_2)(G_K)$. Galois theory reduces this to finding the finite extension $L$ that corresponds to the subgroup $G^{\ell}$ (resp. $G^4$). In principle, we can keep building Kummer $\ell$ extensions until we reach $L$. However this method is not effective unless there exists a numerical criterion to ensure that we can reach $L$ in a certain number of steps. Such criterion can be deduced from Theorem \ref{powerful subgp} if we know the rank of $(\rho_1\oplus \rho_2)(G_K)$.

The main result of this section is Theorem \ref{Lie alg thm} which implies that the rank of $\rho(G_K)_{\ell}$ (resp.~$\rho(G_K)_{4}$) is at most $6$. Then by Theorem \ref{powerful subgp} we give a proof of Theorem \ref{Thm2}. In subsection \ref{section-selfdual-Lie-alg}, we introduce selfdual $\ell$-adic Lie algebras and show that we are reduced to finding the maximal rank of selfdual Lie subalgebras of $\mathfrak{sl}_3(\overline{\Q}_{\ell})$. Then in subsection \ref{section-pf-Thm2}, assuming Theorem \ref{Lie alg thm}, we give a proof to Theorem \ref{Thm2}. The proof to Theorem \ref{Lie alg thm} is in subsection \ref{section-pf-Lie-alg-Thm}. 

\subsection{Selfdual Lie algebras}\label{section-selfdual-Lie-alg}
Recall that for every subgroup $\mathcal{G}$ of $\GL_n(\Z_{\ell})$ there is a Lie algebra $\mathfrak{g}$ attached to it. More specifically, we have the logarithm map
\begin{align*}
\log: 1+{\ell} M_n(\Z_{\ell}) & \longrightarrow M_n(\Z_{\ell})\\
1+y & \longmapsto \sum\limits_{m=1}^{\infty}\frac{(-y)^m}{m}.
\end{align*}
According to \cite[Lemma.~31.1]{Schneider-p-adic-Lie-gp}, there exists an open neighborhood $\mathcal{U}$ of $1_{\mathcal{G}}$ in $\mathcal{G}$ such that the map $\log$  sends $\mathcal{U}$ onto an open neighborhood $\mathfrak{V}$ of $0$ in $\mathfrak{g}$, and the exponential map $\exp(x)=\sum\limits_{m=0}^{\infty}\frac{x^m}{m!}$ gives the local inverse to $\log$.

\begin{Def}
	For each integer $r\geq 0$, we define $\Gamma_r=\{y\in \GL_n(\Z_{\ell}): y\in 1+{\ell}^rM_n(\Z_{\ell})\}$ for $r>0$ and $\Gamma_0=\GL_n(\Z_{\ell})$. If $\mathcal{G}$ is a subgroup of $GL_n(\Z_{\ell})$, then we define $\Gamma_r(\mathcal{G})=\Gamma_r\cap \mathcal{G}$.
\end{Def}

\begin{Prop}\label{expprop}
	For every $r\geq 2$, the $\exp$ and $\log$ maps induce the following group morphisms 
	\begin{align*}
	\log^{(r)}: \Gamma_r/\Gamma_{r+1} &\longrightarrow {\ell}^rM_n(\Z_{\ell})/{\ell}^{r+1}M_n(\Z_{\ell})\\
	1+{\ell}^rx&\longmapsto {\ell}^rx\\
	\end{align*}
	\begin{align*}
	\exp^{(r)}: {\ell}^rM_n(\Z_{\ell})/{\ell}^{r+1}M_n(\Z_{\ell})& \longrightarrow \Gamma_r/\Gamma_{r+1}\\
	{\ell}^ry & \longmapsto 1+{\ell}^ry
	\end{align*}
\end{Prop}
\begin{pf}
	It follows from a straight forward calculation.
\end{pf}

\begin{Def}
	For any subgroup $\mathcal{G}$ of $\GL_n(\Z_{\ell})$, we define its $\bmod \ \ell^r$ rank to be the $\F_{\ell}$-dimension of its image of $\log^{(r)}$. 
\end{Def}

\begin{Prop}\label{modl-rank}
	For every Lie subgroup $\mathcal{G}$ of $\GL_n(\Z_{\ell})$, there exist an integer $r_0=r_0(\mathcal{G})$ such that for all $r\geq r_0$, the $mod \ \ell^r$ ranks of $\mathcal{G}$ is a constant. 
\end{Prop}

\begin{pf}
	We know that the $\Z_{\ell}$ rank of $\mathfrak{g}\cap({\ell}^rM_n(\Z_{\ell}))/\mathfrak{g}\cap({\ell}^{r+1}M_n(\Z_{\ell}))$  are equal to $\dim_{\Q_{\ell}}\mathfrak{g}$ for all $r\geq 0$. Thus the conclusion follows from Proposition \ref{expprop} and the fact for sufficient large $r$, the restriction 
	$$
	\log :\Gamma_r(\mathcal{G}) \to \mathfrak{g}\cap \ell^rM_n(\Z_{\ell})
	$$ 
	is bijective.
\end{pf}

In the following, let $\epsilon=1$ when $\ell=2$ and $0$ otherwise.
\begin{Prop}\label{Gamma-ranks}
	If $\mathcal{G}$ is the image of a $G_K$-representation in $GL_n(\Z_{\ell})$, then
	$$
	rank(\Gamma_{1+\epsilon}(\mathcal{G}))\leq \dim_{\Q_{\ell}}\mathfrak{g}.
	$$
\end{Prop}

\begin{pf}
	Without loss of generality, we can assume that $\mathcal{G}=\Gamma_{1+\epsilon}(\mathcal{G})$. Assume $m=rank(\mathcal{G})$, we construct a minimal set of topological group generators inductively as follow:
	\begin{enumerate}
		\item Let $r_1$ be the least positive integer such that $\Gamma_{r_1}(\mathcal{G})$ has nontrivial image in the quotient $\Gamma_{1+\epsilon}/\Gamma_{r_1+1}$. Then let $a_1$ be an arbitrary element in $\mathcal{G}$ who has nontrivial image in $\Gamma_{1+\epsilon}/\Gamma_{r_1+1}$. Denote by $\mathcal{G}_1=\langle a_1\rangle $ be the topological subgroup generated by $a_1$. 
		\item For every $m \geq j\geq1$, let $r_j\geq r_{j-1}\geq \cdots \geq r_1$ be the least integer such that $\mathcal{G}$ and $\mathcal{G}_{j-1}:=\langle a_1, a_2, \cdots , a_{j-1}\rangle$ have different image in $\Gamma_{1+\epsilon}/\Gamma_{r_j+1}$. Note that such $r_j$ exists since otherwise $\mathcal{G}=\mathcal{G}_j$ which is impossible. Thus we know that $\Gamma_{j_1}(\mathcal{G})$ and $\Gamma_{j_1}(\mathcal{G}_{j-1})$ have different images in $\Gamma_{j_1}/\Gamma_{j_1+1}$. Then we can choose $a_j\in \Gamma_{j_1}(\mathcal{G})$ to be an arbitrary element whose image in $\Gamma_{j_1}/\Gamma_{j_1+1}$ is not contained in that of $\Gamma_{j_1}(\mathcal{G}_{j-1})$.
	\end{enumerate}
	
	According to the above process, we can find a minimal set of topological group generators $\{a_1, \cdots ,a_m\}$. Moreover, note that if $g:=1+x\in 1+\ell^rM_n(\Z_{\ell})$, then 
	$$g^{\ell}=(1+x)^{\ell}=1+\ell x+ O(x^2).$$
	Hence $g^{\ell}= 1+\ell x \ (\bmod \ \ell^{r+2})$. This means if $g_1,g_2\in \Gamma_{r}(\mathcal{G})$ with $r\geq 1+\epsilon$ and $g_1\neq g_2 \ (\bmod \ \Gamma_{r+1})$, then $g_1^{\ell}\neq g_2^{\ell}\ (\bmod \ \Gamma_{j_1+2})$.
	
	With all above discussion, one can see that when $r>r_m$ the set $\{a_1^{r-r_1}, a_2^{r-r_2}, \cdots, a_m^{r-r_m} \}$ induces a basis of the image of $\Gamma_{r}(\mathcal{G})$ in the quotient  $\Gamma_{r}/\Gamma_{r+1}$. Hence in this case the $\bmod \ \ell^r$ rank of $\mathcal{G}$ is always no less than $m$. And this finishes the proof. 
\end{pf}

\begin{Cor}\label{Cor-ranks}
	Let $\mathcal{G}_1$ and $\mathcal{G}_2$ be the image of two $G_K$-representations $\varphi_1$ and $\varphi_2$, respectively. Assume their Lie algebras $\mathfrak{g}_1$ and $\mathfrak{g}_2$ have ranks $r_1 $ and $r_2$, respectively. Then $\Gamma_{1+\epsilon}(img(\varphi_1 \oplus\varphi_2))$ has rank at most $r_1+r_2$.
\end{Cor}

\begin{pf}
	This follows from the above proposition and the fact that $\Gamma_{1+\epsilon}(img(\varphi_1\oplus\varphi_2))$ is a subgroup of $\Gamma_{1+\epsilon}(img(\varphi_1)\oplus img(\varphi_2))$.
\end{pf}

Now if we assume that $\mathcal{G}$ is the image of a strictly selfdual representation $\varphi$, then there exists an invertible matrix $P\in \GL_n(\overline{\Q}_{\ell})$ such that $PgP^{-1}=(g^{-1})^t$ for every $g\in \mathcal{G}$. After taking derivative, we have a selfdual condition for Lie algebras
\begin{align}\label{selfdual-Lie-alg}
PxP^{-1}=-x^t \ \ \  \text{for all } x\in \mathfrak{g}.
\end{align}
After base extension $\mathfrak{g}\otimes \overline{\Q}_{\ell}$, we get a Lie subalgebra of $\mathfrak{gl}_n(\overline{\Q}_{\ell})$ which satisfies the same condition. 

\begin{Def}
	A Lie subalgebra of $\mathfrak{gl}_n(\overline{\Q}_{\ell})$ which satisfies the condition (\ref{selfdual-Lie-alg}) is called a \textit{selfdual Lie algebra}.
\end{Def}
Now we are ready to state the main result in this section. Its proof will be postponed until subsection \ref{classification-lie-alg-section}.
\begin{Thm}{\label{Lie alg thm}}
	If $\mathfrak{g}$ is a selfdual Lie subalgebra of $\mathfrak{sl}_3(\overline{\Q}_{\ell})$, then one of the followings is true.
	\begin{enumerate}
		\item[(a)] $\dim \mathfrak{g}=1$.
		\item[(b)] $\dim \mathfrak{g}=2$ and $\mathfrak{g}$ non-abelian.
		\item[(c)] $\mathfrak{g}\simeq \mathfrak{sl}_2(\overline{\Q}_{\ell})$.
	\end{enumerate}
	In particular, $\dim \mathfrak{g}\leq 3$.
\end{Thm}

\begin{Cor}\label{Cor-rank-3}
	If $\varphi$ is a strictly selfdual representation to $\GL_3(\Z_{\ell})$ with image $\mathcal{G}$, then $\Gamma_{1+\epsilon}(\mathcal{G})$ has rank at most $3$.
\end{Cor}
\begin{pf}
	This follows immediately from the above theorem and Proposition \ref{Gamma-ranks}.
\end{pf}

Now let $\rho_1 $ and $\rho_2$ be two representations such that
\begin{enumerate}
	\item they are both strictly selfdual, and 
	\item they both satisfy condition \ref{cong-trivial-assumption}.
\end{enumerate}   
Recall that $\rho=\rho_1\oplus \rho_2$ has image $\mathcal{G}$.

\begin{Cor}\label{Cor-1}
	There is a filtration of subgroups 
	$$
	\Gamma_{1+\epsilon}(\mathcal{G})=G_1'>G_2'>\cdots>G_{4+\epsilon}'=:V
	$$
	where $G_i'/G_{i+1}'$ is an elementary Abelian $\ell$-group of degree at most $\ell^6$ and $V$ is powerful. Moreover, $[\Gamma_{1+\epsilon}(\mathcal{G}):V]\leq \ell^{18+6\epsilon}$.
\end{Cor}
\begin{pf}
	By Corollary \ref{Cor-rank-3} and Corollary \ref{Cor-ranks}, we know that $\Gamma_{1+\epsilon}(\mathcal{G})$ has rank at most $6$. Thus $\lambda(6)=3$, and the filtration follows from Theorem \ref{powerful subgp} and its remark.
\end{pf}

\begin{Thm}\label{Thm-filtration}
	Let $\mathcal{G}$ be as above, then there is a filtration of subgroups 
	\begin{align*}
	\mathcal{G}&=G_0>G_1>G_2>\cdots>G_{4}=V>G_{5}=V^{\ell} & \text{if $\ell\neq 2$} \\
	\mathcal{G}&=G_0>G_1>G_2>\cdots>G_{6}=V>\cdots >G_{8}=(V^{2})^{2} &\text{ if }\ell=2
	\end{align*}

	such that
	\begin{enumerate}
		\item For every $i\geq 1+\epsilon$, the group $G_i/G_{i+1}$ is an elementary $\ell$-group of rank at most $6$ and $V$ is a powerful subgroup of $\mathcal{G}$.
		\item Every element in $V^{\ell}$ $($resp.~$(V^{2})^2$$)$ is an $\ell$th (resp.~$4$th) power of some other element. Thus $\mathcal{G}/V^{\ell}$ $($resp.~$\mathcal{G}/(V^{2})^2$$)$ surjects onto $\mathcal{G}_{\ell}$ $($$resp.~\mathcal{G}_4)$.
		\item In particular, if $\ell=2$ and $G$ is pro-$2$, then $G_0/G_1$ is elementary of rank at most $4$. 
	\end{enumerate} 
\end{Thm}

\begin{pf}
	If $\ell\neq 2$, we let $G_i=G_i'$ in Corollary \ref{Cor-1} for $i>0$, and if $\ell=2$, we let $G_1=\Gamma_1(\mathcal{G})$ and $G_{i}=G_{i-1}'$ for all $i>1$. Then the first conclusion follows from Corollary \ref{Cor-1}. And Corollary \ref{powerfulcor} implies the second conclusion. To show the last statement, note that $GL_3(\F_2)$ has its strictly upper triangular subgroups as one of its 2-Sylow subgroups. Thus if $\rho_i$ ($i=1,2$) satisfies condition (\ref{cong-trivial-assumption}), its image in $\Gamma_0/\Gamma_1$ is isomorphic to one of the five possible cases:
	$$
	\{1\}, \mathcal{C}_2, \mathcal{C}_2\times \mathcal{C}_2, \mathcal{C}_4, \mathcal{D}_4
	$$ 
	where $\mathcal{C}_n$ is the cyclic group of degree $n$ and $\mathcal{D}_4$ is the dihedral group of degree $8$.
	
	By the next lemma, we know that up to conjugation by an element in $GL_3(\frac{1}{2}\Z_2)$ that the residue image of  ${\rho_i(G_K)}$ is isomorphic to one of the first three cases, hence the last statement follows.
\end{pf}

\begin{Lem}
	Let 
	$$
	P=\left(\begin{smallmatrix} 
	0 & 0 & \frac{1}{2} \\ 
	1 & \frac{1}{2} & \frac{1}{2} \\ 
	1 & 1 & 0
	\end{smallmatrix}  \right) 
	$$
	then 
	$$
	P\left( \begin{smallmatrix} 
	1+2a_{11} & 2a_{12} & 2a_{13} \\ 
	2a_{21} & 1+2a_{22} & 2a_{23} \\ 
	2a_{31} & 2a_{32} & 1+2a_{33}
	\end{smallmatrix} \right) P^{-1}\in 
	\left( \begin{smallmatrix} 
	1 & 0 & a_{31} \\ 
	0 & 1 & a_{21}+a_{31} \\ 
	0 & 0 & 1
	\end{smallmatrix} \right) + 2M_3(\Z_2)
	$$

	$$
	P\left( \begin{smallmatrix} 
	1+2a_{11} & 1+2a_{12} & 2a_{13} \\ 
	2a_{21} & 1+2a_{22} & 2a_{23} \\ 
	2a_{31} & 2a_{32} & 1+2a_{33}
	\end{smallmatrix} \right) P^{-1}\in 
	\left( \begin{smallmatrix} 
	1 & 0 & a_{31} \\ 
	0 & 1 & a_{21}+a_{31} \\ 
	0 & 0 & 1
	\end{smallmatrix} \right) + 2M_3(\Z_2)
	$$
	
	$$
	P\left( \begin{smallmatrix} 
	1+2a_{11} & 1+2a_{12} & 1+2a_{13} \\ 
	2a_{21} & 1+2a_{22} & 1+2a_{23} \\ 
	4a_{31} & 2a_{32} & 1+2a_{33}
	\end{smallmatrix} \right) P^{-1}\in 
	\left( \begin{smallmatrix} 
	1 & 0 & 0 \\ 
	0 & 1 & a_{21} \\ 
	0 & 0 & 1
	\end{smallmatrix} \right) + 2M_3(\Z_2)
	$$

\end{Lem}

\begin{pf}
	Just a calculation.
\end{pf}

\subsection{Proof of Theorem \ref{Thm2}}\label{section-pf-Thm2}

\begin{pf}[Proof of Theorem \ref{Thm2}]
	Now suppose we have two selfdual $G_K$-representations, $\rho_1$ and $\rho_2$, such that $\rho_i^*(2m)\simeq \rho_i$ ($i=1,2$). Then we have 
	$$
	\rho_1\simeq \rho_2\Leftrightarrow \rho_1(-m)\simeq \rho_2(-m).
	$$
	Thus to comparing $\rho_1$ with $\rho_2$, we are reduced to comparing $\rho_1(-m)$ with $\rho_2(-m)$.  Since $(\rho_i(-m))^*=\rho_i^*(m)\simeq \rho_i(-m)$ is strictly selfdual and satisfies condition \ref{cong-trivial-assumption}, we can use Theorem \ref{Thm-filtration} to construct a filtration of subgroups of $\mathcal{G}$. On the other hand, if we build $K_i$ as described above in Theorem \ref{Thm2}, then we have $G_K> \Gal(\overline{K}/K_i)> G_i$ for all $i\geq 1$. Thus there are surjections 
	$$
	\Gal(K_S/K)\twoheadrightarrow \mathcal{G}/G_{5+3\epsilon}\twoheadrightarrow \mathcal{G}_{\ell} (\text{resp}.~\mathcal{G}_4).
	$$
	Now if $T$ is a covering set of $\Gal(K_S/K)$, it is also a covering set of $\mathcal{G}_{\ell}$ $(\text{resp}.~ \mathcal{G}_4)$, then Theorem \ref{Thm2} follows from Proposition \ref{Prop-compare-repns} and Remark \ref{Rmk-covering-set}.
\end{pf}

\subsection{Proof of Theorem \ref{Lie alg thm}}\label{section-pf-Lie-alg-Thm}
\label{classification-lie-alg-section}

In this section, we classify all the selfdual Lie subalgebras of $\mathfrak{sl}_3(\overline{\Q}_{\ell})$ up to conjugacy, and thus prove Theorem \ref{Lie alg thm}. Our arguments are based on detailed calculations. For people who want to skip the calculation details, we provide an outline of our discussion: 
\begin{enumerate}
	\item From Lemma \ref{Lemma-D} to Lemma \ref{delay-proof-1}, we list some basic observations.
	\item After that, we prove that every $2$-dimensional Lie subalgebra of $\mathfrak{sl}_3$ is non-abelian at Proposition \ref{dim2-must-be-nonabelian}, hence prove part $(b)$ of Theorem \ref{Lie alg thm}. 
	\item At Propositions \ref{Prop-Lie-alg-max-dim-3} we show that the selfdual Lie subalgebras of $\mathfrak{sl}_3$ have dimension at most $3$, and all $3$-dimension Lie subalgebras are isomorphic to $\mathfrak{sl}_2$ abstractly. This is the main result of this section since it finishes the proof and Theorem \ref{Lie alg thm}. In order to prove Proposition \ref{Prop-Lie-alg-max-dim-3}
	\begin{enumerate}
		\item From Proposition \ref{Prop-5} to Proposition \ref{Solvable-dim2}, we discuss solvable case. 
		\item At Proposition \ref{Prop-1}, we discuss non-solvable case. 
	\end{enumerate}
\end{enumerate}

Note that if $\mathfrak{g}$ is a selfdual Lie subalgebra of $\mathfrak{gl}_3(\overline{\Q}_{\ell})$, then by condition (\ref{selfdual-Lie-alg})
$$
tr(x)=tr(-x^t)=-tr(x)\Rightarrow tr(x)=0 \text{ for all }x\in \mathfrak{g}
$$ 
thus $\mathfrak{g}$ is a Lie subalgebras of $\mathfrak{sl}_3(\overline{\Q}_{\ell})$.  To simplify the notation, for the rest of this subsection, we simply write $\mathfrak{sl}_3$ for $\mathfrak{sl}_3(\overline{\Q}_{\ell})$. Similarly, we write $\mathfrak{t}$ for $\mathfrak{t}(3,\overline{Q}_{\ell})$, which is the Lie algebra consisting of upper triangular matrices, and  $\mathfrak{n}$ for $\mathfrak{n}(3,\overline{\Q}_{\ell})$, which is the subalgebra of $\mathfrak{sl}_3$ consisting of all strict upper triangular matrices. In this subsection we always assume $\mathfrak{g}$ to be selfdual. In addition, we assume that there is an invertible matrix $P=(p_{ij})\in GL_3(\overline{\Q}_{\ell})$ such that $Px+x^tP=0$ for all $x\in \mathfrak{g}$. Also, we will use the following notation in  \cite{Winternitz-Lie-alg}:

$$
K_1=\frac{1}{2}\left( \begin{smallmatrix} 
1 &  &  \\ 
& -1 &  \\ 
&  &  0
\end{smallmatrix} \right) ,
K_2=\frac{1}{2}\left( \begin{smallmatrix} 
0 & 1 & 0  \\ 
1 & 0 & 0 \\ 
0 & 0 & 0
\end{smallmatrix} \right) ,
K_3=\frac{1}{2}\left( \begin{smallmatrix} 
0 & -1 & 0 \\ 
1 & 0 &  0\\ 
0 & 0 & 0
\end{smallmatrix} \right) , 
D=\left( \begin{smallmatrix} 
1 &  &  \\ 
& 1 &  \\ 
&  & -2
\end{smallmatrix} \right) ,
$$
$$
P_1=\left( \begin{smallmatrix} 
0 & 0 & 1 \\ 
0 & 0 & 0 \\ 
0 & 0 & 0
\end{smallmatrix} \right) ,
P_2=\left( \begin{smallmatrix} 
0 & 0 & 0  \\ 
0 & 0 & 1 \\ 
0 & 0 & 0
\end{smallmatrix} \right) 
R_1=\left( \begin{smallmatrix} 
0 & 0 & 0 \\ 
0 & 0 & 0 \\ 
1 & 0 & 0
\end{smallmatrix} \right) ,
R_2=\left( \begin{smallmatrix} 
0 & 0 & 0 \\ 
0 & 0 & 0 \\ 
0 & 1 & 0
\end{smallmatrix} \right). 
$$ 

Before our dimensional arguments, we have some basic observations. 

\begin{Lem}\label{exact-one-zero}\label{Lemma-D}
	If there is a nonzero diagonal element $\left( \begin{smallmatrix} 
	a &  &  \\ 
	& b &  \\ 
	&  & -a-b
	\end{smallmatrix} \right)$ in $ \mathfrak{g}$, then exactly one of $a$ or $b$ or $a+b$ is $0$. In particular, $D\notin \mathfrak{g}$.
\end{Lem}
\begin{pf}
	To show this, let $A$ to be the given diagonal matrix.  Then,
	$$
	0=PA+A^tP=\left(\begin{smallmatrix} 
	2ap_{11} & (a+b)p_{12} & -bp_{13} \\ 
	(a+b)p_{21} & 2bp_{22} & -ap_{23} \\ 
	-bp_{31} & -ap_{32} & -2(a+b)p_{33}
	\end{smallmatrix}  \right),
	$$
	which proves the lemma by the fact that $P$ is invertible.
\end{pf}

With the same idea and similar calculations, we also have the following lemmas. 

\begin{Lem}\label{2-choose-1}
	At most one of $K_2-K_3$ and $R_2$ is in $\mathfrak{g}$. 
\end{Lem}
\begin{Lem}\label{3-choose-1}\label{delay-proof-1}
	$\mathfrak{g}\cap \mathfrak{n}$ has dimension at most $1$. In particular, at most one of $K_2-K_3$ or $ P_1$ or $P_2$ is in $\mathfrak{g}$.
\end{Lem}

\begin{Rmk}
	By symmetry, $\mathfrak{g}\cap \mathfrak{n}^t$ also has dimension $\leq 1$ also.
\end{Rmk}

Now we discuss dimension. If $\dim \mathfrak{g}=1$, it is trivial, so we start from the $2$-dimensional case and give a proof to part (b) of Theorem \ref{Lie alg thm}.
\begin{Prop}\label{dim2-must-be-nonabelian}
	If $\mathfrak{g}$ is a $2$-dimensional selfdual Lie subalgebra of $\mathfrak{sl}_3$, then $\mathfrak{g}$ is non-abelian. 
\end{Prop}

\begin{pf}
	Suppose $\mathfrak{g}$ is abelian, and let $\{x,y\}$ be its basis. Without loss of generality, we assume that $x$ is of its Jordan form.
	
	If $x$ is diagonalizable, by Lemma \ref{exact-one-zero}, we can write  $
	x=\left( \begin{smallmatrix} 
	1 &  &  \\ 
	& -1 &  \\ 
	&  & 0
	\end{smallmatrix} \right)
	$
	up to scalar. Then $y$ is not diagonal due for Lemma \ref{Lemma-D}. However, since $[x,y]=0$, $y$ must be diagonal, which is impossible. Thus $x$ is not diagonalizable. 
	
	Next suppose 
	$
	x=\left( \begin{smallmatrix} 
	a & 1 & 0 \\ 
	0 & a & 0 \\ 
	0 & 0 & -2a
	\end{smallmatrix} \right)
	$ with $a\neq 0$. Write $y=(y_{ij})_{1\leq i,j\leq 3}$. Then $[x,y]=0$ implies that $y_{13}=y_{21}=y_{23}=y_{31}=y_{32}=0$, and $y_{22}=y_{11}$. Replacing $y$ by $y-\frac{y_{11}}{a}x$ if necessary, we may assume that $y=\left( \begin{smallmatrix} 
	0 & 1 & 0 \\ 
	0 & 0 & 0 \\ 
	0 & 0 & 0
	\end{smallmatrix} \right) $. But then we see that $x-y$ is diagonal, a contradiction. 
	
	If 
	$
	x=\left( \begin{smallmatrix} 
	0 & 1 & 0 \\ 
	0 & 0 & 0  \\ 
	0 & 0 & 0
	\end{smallmatrix} \right),
	$ then $[x,y]=0$ implies $y_{21}=y_{23}=y_{31}=0$, and $y_{22}=y_{11}$. So we can assume $y=\left(\begin{smallmatrix} 
	y_{11} & 0 & y_{13} \\ 
	0 & y_{11} & 0 \\ 
	0 & y_{32} & -2y_{11}
	\end{smallmatrix}  \right)$. Note that if $y_{11}\neq 0$, we may assume that $y_{11}=1$. Then taking $b=0$ if $y_{13}y_{32}=0$ and $-\frac{y_{13}y_{32}}{3}$ otherwise, one can verify that $y+bx$ is diagonalizable, contradiction. Now if $y_{11}=0$, then $y=\left(\begin{smallmatrix} 
	0 & 0 & y_{13} \\ 
	0 & 0 & 0 \\ 
	0 & y_{32} & 0
	\end{smallmatrix}  \right)$,  by Lemma \ref{2-choose-1} and Lemma \ref{3-choose-1}, we have $y_{13}y_{32}\neq 0$. But then we cannot find an invertible matrix $P$ satisfying condition \ref{selfdual-Lie-alg}.
	
	Finally, if 
	$
	x=\left( \begin{smallmatrix}
	0 & 1 & 0 \\ 
	0 & 0 & 1 \\ 
	0 & 0 & 0
	\end{smallmatrix} \right),
	$ we can write $y=\left(\begin{smallmatrix}
	y_{11} & 0 & y_{13} \\ 
	y_{21} & y_{22} & y_{23} \\ 
	y_{31} & y_{32} & -y_{11}-y_{22}
	\end{smallmatrix}  \right)$. Then $[x,y]=0$ implies $y_{21}=y_{23}=y_{31}=y_{32}=0$ and $y_{11}=y_{22}=y_{33}=0$. So we have $y=\left( \begin{smallmatrix} 
	0 & 0 & 1 \\ 
	0 & 0 & 0 \\ 
	0 & 0 & 0
	\end{smallmatrix} \right)$, but then $x,y\in \mathfrak{g}\cap \mathfrak{n}$, contradict Lemma \ref{3-choose-1}. This completes the proof that there is no $2$-dimensional abelian selfdual Lie subalgebra of $\mathfrak{sl}_3$. 
\end{pf}

Suppose the dimension of $\mathfrak{g}$ is at least $3$. We will prove the following proposition which shows part (c) of Theorem \ref{Lie alg thm} and finishes the proof of that theorem. 

\begin{Prop}\label{Prop-Lie-alg-max-dim-3}
	Every selfdual Lie subalgebra $\mathfrak{g}$ of $\mathfrak{sl}_3$ has dimension $\leq 3$. And if  $\dim \mathfrak{g}=3$, then $\mathfrak{g}\simeq \mathfrak{sl}_2$ as abstract Lie algebra.
\end{Prop} 

To prove Proposition \ref{Prop-Lie-alg-max-dim-3}, we will discuss the solvable and non-solvable cases separately (for the related background please see Appendix \ref{Appendix-solvable} or \cite{Humphreys-Lie-Alg}). And Proposition follows immediately from Proposition \ref{Solvable-dim2} and Proposition \ref{Prop-1}.

First, we discuss the solvable cases. Thanks to Lie's theorem (cf. Prop.~\ref{Lie-Thm}), we can assume that $\mathfrak{g}$ is a subalgebra of $\mathfrak{t}$. 
\begin{Prop}\label{Prop-5}
	If $\mathfrak{g}$ is selfdual and solvable, then $\dim \mathfrak{g}\leq 3$.
\end{Prop}

\begin{pf}
	Suppose not, we have $\dim \mathfrak{g}\geq 4$. But then by Lemma \ref{Lemma-D} $\mathfrak{g}\cap \mathfrak{n}$ has dimension at least $2$, which contradicts Lemma \ref{delay-proof-1}. 
\end{pf}

\begin{Prop}\label{Prop-4}
	If $\mathfrak{g}$ is solvable, then $[\mathfrak{g},\mathfrak{g}]=0$ or $\dim[\mathfrak{g},\mathfrak{g}]=1$.
\end{Prop}
\begin{pf}
	This follows from the fact that $[\mathfrak{t},\mathfrak{t}]\subset \mathfrak{n}$ and Lemma \ref{delay-proof-1}.
\end{pf}

\begin{Prop}\label{Solvable-dim2}
	If $\mathfrak{g}$ is a solvable and selfdual Lie subalgebra of $\mathfrak{sl}_3$ then $\dim \mathfrak{g}\leq 2$.
\end{Prop}

\begin{pf}
	Suppose not, we assume $\dim \mathfrak{g}=3$. By Proposition \ref{Prop-4}, we have two possible situations. If $[\mathfrak{g},\mathfrak{g}]=0$, i.e. $\mathfrak{g}$ is Abelian, then according to Lemma \ref{delay-proof-1}, up to scalar there is a unique nonzero element $v$ in $\mathfrak{g}\cap \mathfrak{n}$. Now let $\{v,u_1,u_2\}$ be a fixed basis of $\mathfrak{g}$ as a linear space. Write $v=\left(\begin{smallmatrix} 
	0 & v_{12} & v_{13} \\ 
	0 & 0 & v_{23} \\ 
	0 & 0  & 0
	\end{smallmatrix}  \right)$ and let $u=\left( \begin{smallmatrix} 
	u_{11} & u_{12} & u_{13} \\ 
	0 & u_{22} & u_{23} \\ 
	0 & 0  & -u_{11}-u_{22}
	\end{smallmatrix} \right)$ be an arbitrary element in $\mathfrak{g}\setminus \langle  v \rangle$, then

		$$
		[u,v]=0\Rightarrow	
		\begin{cases}
		u_{11}v_{12}=u_{22}v_{12}\\
		2u_{11}v_{13}+u_{12}v_{23}+u_{22}v_{13}=u_{23}v_{12}\\
		u_{11}v_{23}=-2u_{22}v_{23}.
		\end{cases}
		$$

	If either $v_{12}$ or $v_{23}$ is nonzero, we see that the entries on the main diagonal of $u$ are not independent, and thus there is a nontrivial linear combination of $u_1, u_2$ lying in $\mathfrak{n}$, which implies that $ \dim \mathfrak{g}\cap \mathfrak{n}=2$, contradiction. So $v_{12}=v_{23}=0$. But then $v_{13}\neq 0$, and with the same argument, we still have $\dim \mathfrak{g}\cap \mathfrak{n}=2$. Hence $[\mathfrak{g},\mathfrak{g}]=0$ is impossible. 
	
	Now suppose $\mathfrak{g}$ is non-abelian. We can find linearly independent elements $x,y,v\in \mathfrak{g}$, such that  $[\mathfrak{g},\mathfrak{g}]=\langle v\rangle\subset \mathfrak{n}$. Moreover, we assume at least one of $x$ and $y$ is not commutative with $v$ since otherwise the same arguments as above will imply contradiction. Without loss of generality we let $[v,x]=v$. Then if $[v,y]=bv$ and $[x,y]=cv$, then by replacing $y$ by $y-bx-cv$, we can assume $[x,y]=[v,y]=0$. 
	
	If $v^2\neq 0$, then up to a conjugation by an upper triangular element, we assume that $v=\left(\begin{smallmatrix} 
	0 & 1 & 0 \\ 
	0 & 0 & 1 \\ 
	0 & 0 & 0
	\end{smallmatrix}  \right)$, then the conditions $[v,x]=v$ and $ [v,y]=0$ imply 
	$$
	x=\left( \begin{smallmatrix} 
	-1 & x_{12} & x_{13} \\ 
	0 & 0 & x_{12} \\ 
	0 & 0 & 1
	\end{smallmatrix} \right) \text{ and } y=\left( \begin{smallmatrix} 
	0 & y_{12} & y_{13} \\ 
	0 & 0 & y_{12} \\ 
	0 & 0 & 0
	\end{smallmatrix} \right).$$
	But then $[x,y]=0$ will force $y=0$, contradiction. 
	
	Thus $v^2$ has to be $0$. Up to scaling and conjugation by an upper triangular element, $v$ is one of the following three matrices
	$$
	\left( \begin{smallmatrix} 
	0 & 1 & 0 \\ 
	0 & 0 & 0 \\ 
	0 & 0 & 0
	\end{smallmatrix} \right) ,
	\left( \begin{smallmatrix} 
	0 & 0 & 1 \\ 
	0 & 0 & 0 \\ 
	0 & 0 & 0
	\end{smallmatrix} \right) ,
	\left( \begin{smallmatrix} 
	0 & 0 & 0 \\ 
	0 & 0 & 1 \\ 
	0 & 0 & 0
	\end{smallmatrix} \right).
	$$
	We assume the first. Then by $[v,x]=v$ and $[v,y]=0$, we get 
	$$
	x=\left( \begin{smallmatrix} 
	x_{11} & x_{12} & x_{13} \\ 
	0 & x_{11}+1 & 0 \\ 
	0 & 0 & -2x_{11}-1
	\end{smallmatrix} \right), y=\left( \begin{smallmatrix} 
	y_{11} & y_{12} & y_{13} \\ 
	0 & y_{11} & 0 \\ 
	0 & 0 & -2y_{11}
	\end{smallmatrix} \right).
	$$ 
	By Lemma \ref{delay-proof-1}, $y_{11}\neq 0$, replace $y$ by $\frac{y}{y_{11}}$ and  $x$ by $x-\frac{x_{11}}{y_{11}}y$, we still have $[v,x]=v,[v,y]=[x,y]=0$. Then we have $$
	x=\left( \begin{smallmatrix} 
	0 & x_{12} & x_{13} \\ 
	0 & 1 & 0 \\ 
	0 & 0 &-1
	\end{smallmatrix} \right),
	y=\left( \begin{smallmatrix} 
	1 & y_{12} & y_{13} \\ 
	0 & 1 & 0 \\ 
	0 & 0 & -2
	\end{smallmatrix} \right).$$
	Then $[x,y]=0$ will imply that $y_{13}=3x_{13}$ and $  y_{12}=0$, so we have $y=\left( \begin{smallmatrix} 
	1 & 0 & 3x_{13} \\ 
	0 & 1 & 0 \\ 
	0 & 0 & -2
	\end{smallmatrix} \right)$. But a quick calculation tells us that in this case the invertible $P$ does not exist. For the remaining two choices of $v$, similar arguments give us the same conclusion and finish the proof.
\end{pf}

Second, we talk about non-solvable case.
\begin{Prop}\label{Prop-1}
	If $\mathfrak{g}$ is non-solvable and selfdual, then $\dim \mathfrak{g}\leq 3$. If $\dim \mathfrak{g}= 3$, then $\mathfrak{g}\simeq \mathfrak{sl}_2$.
\end{Prop}
\begin{pf}
	First, we show that $\mathfrak{g}$ has dimension at most $5$. If not, then by linear combination, we can find a Lie subalgebra of $\mathfrak{g}\cap \mathfrak{t}$ which has dimension at least $3$, but this contradicts Theorem \ref{Solvable-dim2}. 
	
	If $\dim \mathfrak{g}\geq 4$, $\mathfrak{g}/Rad(\mathfrak{g})$ is semisimple and nontrivial. By Proposition \ref{ss-alg-dim5} we know the quotient has dimension $3$ and is isomorphic to $\mathfrak{sl}_2$, and thus $\dim Rad(\mathfrak{g})\geq 1$. If we fix $s$ to be an arbitrary nonzero element in $Rad(\mathfrak{g})$, then the following linear map $u\mapsto [s,u]$ from $\mathfrak{g}$ to $ Rad(\mathfrak{g})$  has nonzero kernel. Take $x\neq 0$ to be in that kernel, then $s$ and $x$ spans an abelian dimensional $2$ selfdual Lie subalgebra of $g$. But this contradicts Proposition \ref{dim2-must-be-nonabelian}. Hence, we know that $Rad(\mathfrak{g})=0$ and thus $\mathfrak{g}\simeq \mathfrak{sl}_2$.
\end{pf}

\section{Proof of Theorem \ref{Thm1}}\label{section-pf-Thm1}
In this section, we prove our main theorem \ref{Thm1}. To do this, in subsection \ref{Descend} we show that the symmetric square of the Tate module of elliptic curve $E$ in Theorem \ref{Thm1} can be descended as a $G_{\Q(\sqrt{-3})}$-representation and we give an explicit formula to compute its trace at Frobenius. Then in subsection \ref{Find-covering-set}, by Theorem \ref{Thm2} and the effective Chebotarev density theorem under the Extended RiemannHypothesis (ERH) in \cite{BS-Chebotarev} we find a covering set $T$ by bounding the norm of Frobenius. But that bound is too large for computers to verify Theorem \ref{Thm1}. To find a better $T$, in subsection \ref{Appendix-Burnside-gp} we note that the Galois group $(G_{\Q(\sqrt{-3}),\{2,\infty\}})_{4}\simeq B(2,4)$ (cf. example \ref{eg-burnside}). By studying the conjugacy classes of $B(2,4)$, we prove Theorem \ref{Thm3-trick}, which is the main theorem of this chapter. As a consequence of this theorem, we finally find a set $T$ consisting of no more than $75$ Frobenius as a covering set for Theorem \ref{Thm1}. With the final version of covering set $T$, we are able to finish the calculation in about two weeks. 

Through out this section, we will fix the notation $K$ to be the number field $\Q(\sqrt{-3})$, fix $\mathcal{S}$ and $E$ to be the elliptic surfaces and the elliptic curve involved in Theorem \ref{Thm1}. Also, let $\rho_1$ to be the representation $V_{\ell}$ in Theorem \ref{Thm1}, and take $\rho_2$ to be the descended symmetric square of the Tate module $T_{\ell}(E)$, and $\rho$ to be the direct sum $\rho_1\oplus \rho_2$.

\subsection{Descent of the symmetric square of Tate module}\label{Descend}
Recall that for a number field $F$, and one of its Galois extensions $M/F$, an elliptic curve defined over $M$ is called an $F$-curve if it is isogenous to all its $\Gal(M/F)$ conjugates. Let $\xi$ be a primitive $8$th root of unity and let $b=\left(-\frac{\sqrt{-1}}{4}+\frac{\sqrt{-3}}{8}-\frac{1}{8}\right)$, one can check that the $2$-isogeny
$$
(x,y)\longmapsto \left(\frac{y^2}{2\xi^2x^2}, \frac{y^2(b-x^2)}{2\sqrt{2}\xi^3x^2}\right)
$$
is a morphism from the elliptic curve $E$ in Theorem \ref{Thm1} to its $G_K$-conjugate 
$$
\tilde{E}:Y^2=X^3+(-\sqrt{-1}-1)X^2+\left(\frac{\sqrt{-1}}{4}+\frac{\sqrt{-3}}{8}-\frac{1}{8}\right)X.
$$
Hence $E$ is a $K$-curve. Moreover, one can verify that $E$ does not have complex multiplication. Ribet \cite[$\S$6,7]{Ribet-AV-and-MF} constructs the descent Tate module for every $F$-curve without complex multiplication. We apply Ribet's techniques to our case. 

Suppose $E'$ another elliptic curve and let $\mu:E' \to E$ be an isogeny with dual $\mu^{\vee}$. We write $\mu^{-1}$ to be $\frac{1}{\deg \mu}\mu^{\vee}\in Hom(E,E')\otimes \Q$. For every element $\sigma\in G_K$, we denote by $E^{\sigma}$ the conjugation of $E$ by $\sigma$, and fix an isogeny $\mu_{\sigma}:E^{\sigma}\to E$. Then the following map
\begin{align*}
c: G_{K}\times G_{K} & \longrightarrow \Q\\
(\sigma, \tau) & \mapsto \mu_{\sigma}\mu_{\tau}^{\sigma}\mu_{\sigma\tau}^{-1}
\end{align*}
is a well-defined (recall that $E$ does not have complex multiplication) $2$-cocycle. By the following Proposition \ref{Prop-Ribet-Tate} we know that $c$ is a $2$-boundary, i.e. there exists $\alpha:G_K\to \overline{\Q}^*$, such that 
$$
c(\sigma, \tau)=\frac{\alpha(\sigma)\alpha(\tau)}{\alpha(\sigma\tau)}.
$$
\begin{Prop}\cite[Thm.~6.3]{Ribet-AV-and-MF}\label{Prop-Ribet-Tate}
	$H^2(G_K, \overline{\Q}^*)=0$, where $\overline{\Q}^*$ has the trivial $G_K$-action.
\end{Prop}

Now we define a $G_{K}$-action $\phi$ on $\overline{\Q}_l\otimes T_{{E}}$ by 
$$
\phi(\sigma)(1\otimes x)=\alpha^{-1}(\sigma)\otimes \mu_{\sigma}(x^{\sigma}).
$$
It is a well-defined for $
\phi(\sigma)\phi(\tau)\phi(\sigma\tau)^{-1}=c(\sigma,\tau)(\frac{\alpha(\sigma)\alpha(\tau)}{\alpha(\sigma\tau)})^{-1}=1.
$

Every conjugation $E^{\sigma}$ is either isomorphic to $E$ or $\tilde{E}$ over $K$. Let $\mu:\tilde{E}\to E$ be a $2$-isogeny and take 

	$$
	\mu_{\sigma}=
	\begin{cases}
	1 \text{ if }\sigma \in G_{K(\sqrt{-1})}\\
	\mu \text{ if } \sigma\notin G_{K(\sqrt{-1})}.
	\end{cases}
	$$

By calculation, we have 
\begin{enumerate}
	\item[(1)] If $\sigma, \tau\in G_{K(\sqrt{-1})}$, then $c(\sigma, \tau)=\mu_{\sigma}\mu_{\tau}^{\sigma}\mu_{\sigma\tau}^{-1}=1\circ1\circ 1=1$;
	\item[(2)] If $\sigma\in G_{K(\sqrt{-1})}, \tau\notin G_{K(\sqrt{-1})}$, then $c(\sigma, \tau)=\mu_{\sigma}\mu_{\tau}^{\sigma}\mu_{\sigma\tau}^{-1}=1\circ\mu\circ \mu^{-1}=1$;
	\item[(3)] If $\sigma\notin G_{K(\sqrt{-1})}, \tau\in G_{K(\sqrt{-1})}$, then $c(\sigma, \tau)=\mu_{\sigma}\mu_{\tau}^{\sigma}\mu_{\sigma\tau}^{-1}=\mu\circ 1\circ \mu^{-1}=1$;
	\item[(4)] If $\sigma, \tau\notin G_{K(\sqrt{-1})}$, then $c(\sigma, \tau)=\mu_{\sigma}\mu_{\tau}^{\sigma}\mu_{\sigma\tau}^{-1}=\mu\circ \mu^{\sigma}\circ 1=\pm\deg (\mu)=\pm2$.
\end{enumerate} 
Thus, let 
$$
\alpha(\sigma)=\begin{cases}
1 \text{ if }\sigma\in G_L\\
\sqrt{2} \text{ if } \sigma \notin G_L
\end{cases}
$$
then we can descend $T_{{E}}$ as $G_K$-representation as following:

\begin{align*}
\phi: G_K & \to End(T_{{E}})\otimes \Q_{\ell}\\
\sigma & \mapsto \begin{cases}
x\mapsto x^{\sigma}\ \ &\text{if }\sigma\in G_{K(\sqrt{-1})}\\
x\mapsto \frac{\mu(x^{\sigma})}{\sqrt{2}} \ \ &\text{if }\sigma\notin G_{K(\sqrt{-1})}
\end{cases}.
\end{align*}

\begin{Rmk}
	The image of $\phi$ is in $GL_2(\frac{1}{\sqrt{2}}\Z_{\ell})$. But since we only care about $\rho_2$ which is the symmetric square of $\phi$, we know that the image of $\rho_2$ is still in $GL_3(\Z_{\ell})$. By the same reason, $\rho_2$ will not change if we choose $\alpha(\sigma)=-\sqrt{2}$ for $\sigma\notin  G_{K(\sqrt{-1})}$.
\end{Rmk}

Now we need to compute the trace of $\rho_2=Sym^2(\phi)$ for every Frobenius $F_{\mathfrak{p}}$ over $K$. If we denote by $\pi$ the representation induced by $T_{{E}}$, it's easy to see that when $\mathfrak{p}$ splits in $K(\sqrt{-1})$, we have $
tr(Sym^2(\phi)(F_{\mathfrak{p}}))=tr(Sym^2(\pi)(F_{\mathfrak{p}})).$
Now assume $\mathfrak{p}$ is inert in $K(\sqrt{-1})$ , with $\mathfrak{P}$ lying above it. With a proper choice of the basis of $T_{{E}}\otimes \overline{\Q}_{\ell}$ such that $\phi(F_{\mathfrak{p}})=\left(\begin{smallmatrix}
a & * \\ 
0 & b
\end{smallmatrix}  \right) $, we get 
$$
tr(Sym^2(\phi)(F_{\mathfrak{p}}))=a^2+b^2+ab=tr(\pi(F_{\mathfrak{P}}))+ab.
$$ Since $(ab)^2=\det(\pi(F_{\mathfrak{P}}))=N_{K/\Q}(\mathfrak{p})^2$, we know that $\det(\phi(F_{\mathfrak{p}}))=ab=\pm N_{K/\Q}({\mathfrak{p}})$. To determine the sign, we use the idea in the proof of \cite[chap.~\RNum{5}, Prop.~2.3, Thm.~2.4]{Silverman1}: if we consider the reduced curve of $E$ on the residue field at $\mathfrak{p}$, then the determinant of $\sqrt{2} \phi(F_{\mathfrak{p}})$ can be explained as the degree of isogeny $\mu \circ F_{\mathfrak{p}}$, and hence it must be positive. As conclusion, we have 

\begin{align}\label{trace-formula}
tr(Sym^2(\phi)(F_{\mathfrak{p}}))=
\begin{cases}
tr(\pi(F_{\mathfrak{p}}))^2-N_{K/\Q}(\mathfrak{p}) \text{ when } \mathfrak{p} \text{ splits}\\
tr(\pi(F_{\mathfrak{P}}))+N_{K/\Q}(\mathfrak{p})	\text{ otherwise }
\end{cases}
\end{align}

In particular, since $(0,0)$ is a $2$-torsion point of $E$, we know that $tr(\pi(F_{\mathfrak{P}}))=0\ (\bmod 2)$ for all  finite prime ideals $\mathfrak{P}$ in ${K(\sqrt{-1})}$. Hence we know that all $\bmod 2$ characteristic polynomials of $Sym^2(\phi)$ are equal to
\begin{equation}\label{eqn-2}
t^3+t^2+t+1=(t-1)^3 \quad (\bmod\ 2),
\end{equation}
thus $\rho_2(-1)$ satisfies the condition (\ref{cong-trivial-assumption}).

\subsection{Finding a covering set by Theorem \ref{Thm2}}\label{Find-covering-set}
When $\ell=2$, by (\ref{eqn-1}) and (\ref{eqn-2}) we know that $\rho_1(-1)$ and $\rho_2(-1)$ are congruent trivial (condition (\ref{cong-trivial-assumption})). Moreover, since both $\mathcal{S}$ and $E$ are smooth outside primes above $2$. The Galois representations $\rho_1$ and $ \rho_2$ are unramified outside of the finite set $S=\{2,\infty\}$. According to Theorem \ref{Thm2} we can find a covering set by the algorithm below:
\begin{enumerate}
	\item Take $K_0=K$, then for every $i\geq 0$, list all quadratic extensions $L/K_i$ which satisfy the following conditions:
	\begin{enumerate}
		\item Unramified outside $S=\{2,\infty\}$;
		\item For every prime place $\mathfrak{p}$ in $K$ not dividing $2$, and for every prime place $\mathfrak{P}$ in $L$ above $\mathfrak{p}$, the corresponding local field extension has Galois group of exponent no greater than $4$.
	\end{enumerate}  
	\item Take $K_{i+1}$ to be the compost of all the $L$ listed above. And build $K_{i+2}$ inductively, until either
	\begin{enumerate}
		\item There is no such quadratic extension $L/K_i$ satisfying the conditions in step (1), which means this $K_i$ is the maximal exponent $4$ pro-$2$ extension above $K$ which is unramified outside $S$, or;
		\item $i=8$, which means $K_i=K_S$ in Theorem \ref{Thm2}.
	\end{enumerate} 
	\item Denote by $K_{max}$ the field which the above process ends up to, and use the effective Chebotarev density theorem to find a bound $B$ of the absolute norm of Frobenius classes for the covering set of $\Gal(K_{max}/K)$. Then 
	$$T=\{F_{\mathfrak{p}}| N_{K/\Q}(\mathfrak{p})\leq B\}$$
	is sufficient to be a covering set. 
\end{enumerate}

In our situation, we have $K_1=\Q(\zeta_{24})$. Then $K_2/K_1$ has Galois group isomorphic to $(\Z/2\Z)^{\oplus 5}$. But now $[K_2:\Q]=2^8=256$, which means $K_2$ is very hard to be constructed via computers. To get more information without extending $K_2$, recall that $\Gal(K_{max}/K)$ is a quotient of $\Gal(K^{ur}_{\infty,2}(2)/K)_4$ where $K^{ur}_{\infty,2}(2)$ is the maximal pro-$2$ extension above $K$ and unramified outside $S=\{2,\infty\}$. Hence we are reduced to finding a covering set of $\Gal(K^{ur}_{\infty,2}(2)/K)_4$. By \cite[Thm.~2]{Jossey-2-ext-ur-outside-2}, and calculations with help of computers, we find that:
\begin{enumerate}
	\item  $\Gal(K^{ur}_{\infty,2}(2)/K)$ is isomorphic to the free group generated by two elements, hence $\Gal(K^{ur}_{\infty,2}(2)/K)_4\simeq B(2,4)$ (cf. Example \ref{eg-burnside}).
	\item The natural quotient map $\Gal(K^{ur}_{\infty,2}(2)/K)_4\rightarrow \Gal(K_2/K)$ has kernel isomorphic to $(\Z/2)^{\oplus 5}$, so $K_{max}=K_3$, and $\Gal(K_3/K_2)$ is a subgroup of $(\Z/2\Z)^{\oplus 5}$. 
	\item Since $K_3/K_2$ is an abelian extension, by Conductor-Discriminant formula, we get $disc(K_3)\leq 2^{109568}3^{4096} $. Assume Extended Riemann Hypotheses and apply the the following Theorem \ref{Chebotarev} and its remark to find a covering set 
	$$T:=\{\mathfrak{p}\in K| N_{K/\Q}({\mathfrak{p}})\leq 7.03\times 10^9\}.$$
\end{enumerate}

\begin{Thm}\cite[Thm.~5.1]{BS-Chebotarev} \label{Chebotarev}
	Let $\Delta$ be the discriminant of $L/K$ and take $n=[L:K]$. If we assume the Extended Riemann Hypotheses, then 
	$$
	B(L/K)\leq (4\log \Delta+2.5n+5)^2.
	$$ 
\end{Thm}

\begin{Rmk}
	In practice, instead of using the above inequality, one can write a code to find a sharper bound with the idea in \cite[$\S$~3,4]{BS-Chebotarev}. This is what we did in our case. 
\end{Rmk}

\subsection{Conjugacy classes in $B(2,4)$, proof of Theorem \ref{Thm4}}\label{Appendix-Burnside-gp}
The covering set $T=\{\mathfrak{p}|N_{K/\Q}({\mathfrak{p}})\leq 7.03\times 10^9 \}$ is sufficient large to verify Theorem \ref{Thm1} by Faltings-Serre method. But even with a Xeon-E3 CPUs, it would take years to finish calculations. In this subsection, we reduce the size of $T$. The main result of this subsection is Theorem \ref{Thm3-trick}. With this theorem we can find a new covering set consisting only $75$ prime ideals without requiring Extended Riemann Hypotheses. Our method is based on the group structure of $B(2,4)$ and Galois theory. For this subsection, we loose our restriction on $K$ and $\rho_i$ ($i=1,2$) by supposing they are as stated in Theorem \ref{Thm4}. Hence taking $K=\Q(\sqrt{-3})$ and $\rho_1$ and $\rho_2$ to be the selfdual representations induced by $V_{\ell}$ and $Sym^2(T_E)$, we prove Theorem \ref{Thm1}.

Recall that if we have a covering set of $\Gal(K_{2,\infty}^{ur}/K)$, then to verify Theorem \ref{Thm1}, we need to compare the trace of $\rho_1$ and $\rho_2$ for every element in the covering set. But in fact if two elements $x$ and $y$ are conjugate in $\Gal(K^{ur}_{2,\infty})$, then $\rho_i(x)$ and $\rho_i(y)$ have the same trace ($i=1,2$). Hence to verify Theorem \ref{Thm1}, we only need to check if $\rho_1$ and $\rho_2$ share the same trace for every representative of conjugate classes in $T$. Therefore, a covering set of the conjugate classes of $\Gal(K^{ur}_{\infty,2}(2)/K)_4\simeq B(2,4)$ will be sufficient to verify Theorem \ref{Thm1}. 

So we are reduced to find a covering set (still denoted by $T$) of the (conjugate) classes in $B(2,4)$. If we denote by $C_1,\cdots ,C_{88}$ the $88$ classes in $B(2,4)$, then to find a covering set of those $C_i$'s, we need to tell which class $F_{\mathfrak{p}}$ belongs to for every given prime ideal $\mathfrak{p}$ of $\OO_K$. This is not easy since the isomorphism between $\Gal(K^{ur}_{2,\infty})$ and $B(2,4)$ is not explicit. However, suppose we can find $88$ distinct Frobenius elements $F_{\mathfrak{p}_i}$ and show that they belong to distinct classes, then the set $\{\mathfrak{p}_i\}$ is a desired covering set. 

In order to differentiate non-conjugate Frobenius elements, we note that if $N$ is a normal subgroup of $B(2,4)$, then for each $i$ either $C_i\subset N$ or $C_i\cap N=\emptyset$. Thus we can differentiate two classes $C_i$ and $C_j$ if there is a normal subgroup $N$ which contains one of the two classes and is disjoint with the other. Hence, if we can find a finite set $\mathcal{N}$ of normal subgroups of $B(2,4)$ such that every conjugate class has a unique ``intersection pattern'' (Definition \ref{Def-pattern-c}) with respect to elements in $\mathcal{N}$. Then $\mathcal{N}$ can be used to differentiate all the classes. On the other hand, since every Frobenius element $F_{\mathfrak{p}}$ represents a conjugate classes of elements in the Galois group, we can tell that $F_{\mathfrak{p}_1}$ and $F_{\mathfrak{p}_2}$ are not in the same conjugate class if one of them is in a normal subgroup $N$ and another is not. Therefore, if we are able to find $88$ $F_{p}$'s so that they have distinct ``intersection patterns'' with respect to elements in $\mathcal{N}$ above, then these $88$ $F_{\mathfrak{p}}$'s will provide a desired covering set. To realize this method effectively, we need the following notations and definitions. 

\begin{Def}\label{Def-pattern-c}
	For a group $G$, let $\mathcal{C}=\{C_1,\cdots, C_r\}$ be a set of (conjugate) classes and let  $\mathcal{N}=\{N_1,\cdots ,N_s\}$ be an ordered set consisting of normal subgroups of $G$. For each $1\leq i\leq r$ and $1\leq j\leq s$, take 
	$$
	\delta_{i,j}=
	\begin{cases}
	1 \ \ \text{if }C_i\subset N_j\\
	0 \ \ \text{otherwise}
	\end{cases}.
	$$
	Then for a fixed $C_i$, we define the vector $P(C_i,\mathcal{N})=(\delta_{i,1},\delta_{i,2}, \cdots ,\delta_{i,s})$ the \textit{pattern of $C_i$ with respect to $\mathcal{N}$}. If there is no confusion, we will simply call it the \textit{pattern of $C_i$} and write $P(C_i)$ for short.
\end{Def} 

\begin{Def}\label{Def-pattern-p}
	Let $F$ be a number field, denote by $\mathcal{L}=\{L_1,\cdots ,L_s\}$ an ordered set of finite Galois extensions $F$. For a prime ideal $\mathfrak{p}$ in $F$, we define 
	$$
	\delta(\mathfrak{p}, L_j)=
	\begin{cases}
	1\ \ \text{if }\mathfrak{p} \text{ is toltally split in }L_j\\
	0\ \ \text{otherwise}
	\end{cases}.
	$$
	The vector $P(\mathfrak{p}, \mathcal{L})=(\delta(\mathfrak{p}, L_1), \cdots, \delta(\mathfrak{p}, L_s))$  is defined to be the \textit{pattern of $\mathfrak{p}$ with respect to $\mathcal{L}$}. When there is no confusion, we call it the \textit{pattern} of $\mathfrak{p}$ and simply write $P(\mathfrak{p})$.
\end{Def}

We can rephrase our method by the definitions above. We want to find an order set $\mathcal{N}$ of normal subgroups of $B(2,4)\simeq \Gal(K^{ur}_{\infty,2}(2)/K)_4$ such that every class has a unique pattern with respect $\mathcal{N}$. Then by Galois theory, for every $N_j\in \mathcal{N}$, we denote by $L_j/K$ the corresponding Galois intermediate extension of $K^{ur}_{2,\infty}/K$. Then take $\mathcal{L}=\{L_{j}\}$ to be the corresponding ordered set. Then we know that 
$$
F_{\mathfrak{p}}\in C_i \text{ if and only if } P(\mathfrak{p}, \mathcal{L})=P(C_i, \mathcal{N}).
$$
Hence if we want to find a desired covering set, we only need to factor the prime ideals of $\OO_K$ one-by-one until we find $88$ primes with distinct patterns. Based on this idea, we will state the main theorem of this section (Theorem \ref{Thm3-trick}) after some definitions and lemmas. 

\begin{Def}
	Let $F$ be a number field, a Galois extension $L/F$ is called an \textit{exponent $4$ extension} if its Galois group $\Gal(L/F)$ has exponent $4$.
\end{Def}

\begin{Lem}\label{Lem-7-exp4-ext}
	For every number field in Theorem \ref{Thm4}, there are exact $7$ quartic exponent $4$ Galois extensions which are unramified outside $\{2, \infty \}$. In particular, among them, there is exact one biquadratic extension.
\end{Lem}

\begin{pf}
	According to \cite[Theorem~2]{Jossey-2-ext-ur-outside-2} we know that for each $K$ satisfying the condition in Theorem \ref{Thm4}, its maximal pro-$2$ extension which is unramified outside $\{2, \infty\}$ (i.e. $K^{ur}_{\infty,2}(2)/K$) has a free pro-$2$ Galois group generated by $2$ elements. Hence $\Gal(K^{ur}_{\infty,2}(2)/K)_4\simeq B(2,4)$. The Burnside group $B(2,4)$ has order $2^{12}$. One can find by calculation that there are exact $7$ order $2^{10}$ normal subgroups $N$ of $B(2,4)$, and exact one of them satisfies $B(2,4)/N\simeq (\Z/2\Z)^{\oplus 2} $. This finishes the proof. 
\end{pf}

\begin{Thm}\label{Thm3-trick}
	Let $K$ and $\rho_1$, $\rho_2$ be as in Theorem \ref{Thm4}. Denote by $L_i$ $(i=1, \cdots , 7)$ the $7$ Galois quartic extensions of $K$ which are unramified outside $\{2, \infty \}$, and let $L_1$ be the unique biquadratic extension of $K$. For each $1\leq i \leq 7$, let $\mathcal{L}_1=\{M_{i,s} \}$ be an ordered set of all exponent $4$ Galois extensions of $L_i$ such that $M_{i,s}/L_i$ is unramified outside $\{2, \infty \}$, and $[M_{i,s}:L_i]\leq 2^{r_i}$ with $r_i=3$ if $i=1$ and $2$ otherwise. Then we have the followings.	
	
	\begin{enumerate}
		\item There are $63$ prime ideals $\mathfrak{p}_h$ $(h=1, \cdots, 63)$ of $\OO_K$ satisfying the followings. 
		\begin{enumerate}
			\item $\mathfrak{p}_h\nmid 2$ and $\mathfrak{p}_h$ is totally split in $\OO_{L_1}$.
			\item Let $U_1=\{\mathfrak{P}_{h,j}\}$ consist of all prime ideals $\mathfrak{P}_{h,j} $ of $\OO_{L_1}$ such that $\mathfrak{P}_{h,j}$ is lying above $\mathfrak{p}_h$. Then the set of patterns $\{P(\mathfrak{P}_{h,j}), \mathcal{L}_1 \}$ has $204$ elements. 
		\end{enumerate}
		\item For each $i\in \{2, \cdots , 7 \}$, there are $2$ prime ideals $\mathfrak{p}_h$ $(h=1, 2)$ of $\OO_K$ satisfying the followings. 
		\begin{enumerate}
			\item $\mathfrak{p}_h\nmid 2$ and $\mathfrak{p}_h$ is totally split in $\OO_{L_i}$.
			\item Let $U_i=\{\mathfrak{P}_{h,j}\}$ consist of all prime ideals $\mathfrak{P}_{h,j} $ of $\OO_{L_i}$ such that $\mathfrak{P}_{h,j}$ is lying above $\mathfrak{p}_h$. Then the set of patterns $\{P(\mathfrak{P}_{h,j}), \mathcal{L}_i \}$ has $8$ elements. 
		\end{enumerate}		
	\end{enumerate}
	Let $T$ be the collection of all the prime ideals of $\OO_K$ stated in $(1)$ and $(2)$, then 
	$$
	\rho_1\sim\rho_2 \quad \text{if and only if} \quad tr(\rho_1(F_{\mathfrak{p}}))=tr(\rho_2(F_{\mathfrak{p}}))\ \text{for all }\mathfrak{p}\in T.
	$$
	In particular, when $K=\Q(\sqrt{-2})$, the set $T$ is given by Table \ref{List2}. When $K=\Q{\sqrt{-3}}$, $T$ is given by Table \ref{List3}.
\end{Thm}

\begin{Rmk}
	The above results do not change if we change the order of $L_i$ in $\mathcal{L}$.
\end{Rmk}

\begin{pf}[Proof of Theorem \ref{Thm3-trick}]
	Let $\mathcal{N}=\{N_1, \cdots , N_7\}$ be the ordered set consisting of all normal subgroups of order $\geq 2^{10}$ (i.e. those corresponding to Galois extensions over $K$ and having absolute extension degree $\leq 2^3=8$). In particular we fix $N_{1}$  to be the unique group such that $B(2,4)/N_{0}\simeq (\Z/2\Z)^{\oplus 2}$. By calculating the patterns of every $C_i$ with respect to $\mathcal{N}$, we partially differentiate them in to $7$ subsets $TC_1, \cdots , TC_7$. In Table \ref{Table1} we list all the $7$ subsets by writing down their classes the common pattern of the classes in each subset. The characteristic of the sets $TC_i$ ($i=0,\cdots 7$) can be read from this table.

	\begin{table}[h!]
		\begin{center}
			\begin{tabular}{| c | c |   c | c | c | c | c |}
				\hline
				$TC_i$	& classes in the set	& Description of the pattern 	\\ \hline
				
				$TC_1$ 	& $C_1,\cdots, C_{64}$	& $\delta_{i,1}=1$ 	\\ \hline
				
				$TC_2$ & $C_{65}, C_{71}, C_{84}, C_{85}$	& \multirow{6}{*}{\makecell{$\delta_{i,j}=1$ for exact \\ one of $2\leq j\leq 7$}}\\ \cline{1-2}
				
				$TC_3$ & $C_{66}, C_{72}, C_{74}, C_{78}$	&	\\ \cline{1-2}
				
				$TC_4$ & $C_{67}, C_{68}, C_{87}, C_{88}$	&	\\ \cline{1-2}
				
				$TC_5$ & $C_{69}, C_{80}, C_{82}, C_{86}$	&	\\ \cline{1-2}
				
				$TC_6$ & $C_{70}, C_{75}, C_{81}, C_{83}$	&	\\ \cline{1-2}
				
				$TC_7$ & $C_{73}, C_{76}, C_{77}, C_{79}$	&	\\ \hline

			\end{tabular}
			
			\caption{Classes of the same patterns with respect to $\mathcal{N}=\{N_1,\cdots ,N_{7}\}$}
			\label{Table1}

		\end{center}
	\end{table}
	
	\begin{Rmk}
		\begin{enumerate}
			\item[(a)] Theoretically, adding more normal subgroups of $B(2,4)$ to $\mathcal{N}$ will be helpful to refine Table \ref{Table1}. However, it is hard in practice since the normal subgroups we do not consider are corresponding to extensions of degree $\geq 2^7=128$, which take too much time to construct. Hence we finally take our $\mathcal{N}$ to be as above. 
			\item[(b)] Changing the order of elements in $\mathcal{N}$ does not change the classification in Table \ref{Table1}.
		\end{enumerate}
	\end{Rmk}
	
	Now in order to differentiate the classes in $TC_1$, we let $\mathcal{N}_1$ to be the ordered set of all normal subgroups of $N_{1}$ with order $\geq 2^7$. Then applying the same idea as above except with $\mathcal{N}$ replaced by $\mathcal{N}_1$, we find the followings by calculation. 
	\begin{enumerate}
		\item A conjugate class of $B(2,4)$ may not still be a single conjugate class with respect to $N_1$. In fact, for each $C_i$ for $i=1, \cdots , 64$, we have $C_i=\bigsqcup C_{i,k}$ is a disjoint union of several conjugate classes of $N_1$. By calculation, there are totally $208$ sub classes \{$C_{i,k}$\} spitted from $C_1, \cdots , C_{64}$. 
		\item By computing all the patterns $P(C_{i,k}, \mathcal{N}_1)$, we find $204$ distinct patterns. In particular, except that the subclasses from $C_{63}$ have the same patterns with the subclasses from $C_{64}$, each other subclass has a unique pattern with respect to $\mathcal{N}_1$.
		\item If $g\in C_{63}$, then $g^{-1}\in C_{64}$.
	\end{enumerate}
	Therefore, to find a covering set of the set $\{C_{i}\}_{i=1, \cdots, 64}$, we are reduced to finding a covering set of $\{C_{i,k}\}_{i=1, \cdots, 64}$. For the later, take $\mathcal{L}_1$ to be the ordered set in Theorem \ref{Thm3-trick}, then we need to find $204$ prime ideals of $\OO_K$ which totally splits in $L_1$, and have distinct patterns with respect $\mathcal{L}_1$. Note that there is no need to differentiate $C_{63}$ from $C_{64}$ since by part $(3)$ above, once we find an element in either one of the two classes, we automatically find an element in another. This proves part (1) of Theorem \ref{Thm3-trick}.  
	
	To differentiate the classes in each of $TC_i$ $(i=2, \cdots , 7)$, we let $\mathcal{N}_i$ be the ordered set of all normal subgroups of $N_i$ with order $\geq 2^8$. Then the classes in each $T_i$ will split into $16$ sub classes of $N_i$. By computing their patterns with respect to $\mathcal{N}_i$, we have $8$ distinct patterns. Moreover, if $g$ is in a sub class $C$ then $g^{-1}$ is the sub class $C'$ such that $C'\neq C$ but $P(C, \mathcal{N}_i)=P(C', \mathcal{N}_i)$. Hence to find a covering set to each of $T_i$ ($i=2, \cdots , 7$), we take $\mathcal{L}_i$ to be the corresponding ordered set in Theorem \ref{Thm3-trick}, and then find $8$ prime ideals of $\OO_K$ which totally splits in $L_i$ and have distinct patterns with respect to $\mathcal{L}_i$. This proves part (2) of Theorem \ref{Thm3-trick}.
	
	Finally, if we denote by $T_i$ ($i=1, \cdots , 7$), the covering $TC_i$, and take $T$ to be the union of all $T_i$. Then $T$ is sufficiently large in the sense that every conjugate class has a representative in terms of an element in $T$ or the inverse of an element in $T$. As a conclusion, by Faltings-Serre method, to verify whether $\rho_1$ is equivalent to $\rho_2$, we only need to test if they have the same characteristic polynomial for every element in $T$. In particular, when $K=\Q(\sqrt{-2})$ and $K=\sqrt{-3}$, the corresponding sets $T$ are given by Tables \ref{List2} and \ref{List3} respectively. This finishes the proof. 
\end{pf}

\begin{table}[h!]
	\begin{center}
		\begin{tabular}{| c | c | }
			\hline
			$T_i$	& prime integers lying below $\mathfrak{p}$ 	\\ \hline
			
			$T_1$ 	& \makecell{ 439, 503, 607, 823, 1231, 1399, 1423, 3049, 3089\\ 3449, 3823, 3967, 4057, 4177, 4201, 4217,
				4409, 4937\\ 5737, 6121, 6353, 6553, 7793, 9377, 9473, 9769, 11113\\ 11969, 12241, 16433, 18593, 25409, 
				26993, 27809, 67217\\ 67489, 68449, 126641, 132929, 268817, 392737}	\\ \hline
			
			$T_2$ 	& \multirow{6}{*}{\makecell{29,67,97,137,139,193,251,283}}\\ 
			
			$T_3$ 	&	\\ 
			
			$T_4$ 	&	\\ 
			
			$T_5$ 	&	\\ 
			
			$T_6$ 	&	\\ 
			
			$T_7$ 	&	\\ \hline

		\end{tabular}
		
		\caption{The primes in $T$ when $K=\Q(\sqrt{-2})$}
		\label{List2}

	\end{center}
\end{table}

\begin{table}[h!]
	\begin{center}
		\begin{tabular}{| c | c | }
			\hline
			$T_i$	& prime integers lying below $\mathfrak{p}$ 	\\ \hline
			
			$T_1$ 	& \makecell{ 419, 461, 587, 617, 647, 653, 761, 911, 929, 983, 1439, 2273\\ 2521, 3023, 3793, 3889, 4297, 4513, 4969, 5113, 6337, 6673\\ 7393, 8161, 8329, 8353, 8641, 9049, 9337, 9721, 10369\\ 10729, 11113, 11161, 12577, 13873, 14713, 15121, 15913\\ 19777, 21193, 25537, 31393, 40177, 57697, 71233, 74353\\ 87697, 98641, 100801, 104593, 115153, 234721}	\\ \hline
			
			$T_2$ 	& \multirow{6}{*}{\makecell{37,127,181,199,211,271,379,523,619,631}}\\ 
			
			$T_3$ 	&	\\ 
			
			$T_4$ 	&	\\ 
			
			$T_5$ 	&	\\ 
			
			$T_6$ 	&	\\ 
			
			$T_7$ 	&	\\ \hline

		\end{tabular}
		
		\caption{The primes in $T$ when $K=\Q(\sqrt{-3})$}
		\label{List3}

	\end{center}
\end{table}

\begin{pf}[Proof of Theorem \ref{Thm4}]
	It immediately follows from Theorem \ref{Thm3-trick}.
\end{pf}

\begin{pf}[Proof of Theorem \ref{Thm1}]
	First by (\ref{eqn-1}) and (\ref{eqn-2}) we know that $\rho_1(-1)$ and $\rho_2(-1)$ are congruent trivial. Then by running code on twelve CPUs, it takes about $15$ days to finish the comparison. As a result, we know that the two representations in Theorem \ref{Thm1} are equivalent. Then by the fact that $E$ does not have complex multiplication and apply Serre's open image theorem (Theorem \ref{Thm-Serre-open-image}) and its Corollary \ref{Cor-Sym-irr} we know that the two representations in Theorem \ref{Thm1} are both irreducible, thus they are isomorphic to each other. This finishes the proof of Theorem \ref{Thm1}.
\end{pf}

\appendix

\section{Backgrounds of Lie algebras}\label{Appendix-Lie-alg}
This section is devoted to fill the backgrounds of Lie algebras that are needed in our paper.

\subsection{Solvable Lie algebras}\label{Appendix-solvable}
In this subsection, $\mathfrak{g}$ is a Lie algebra which is not necessary selfdual. 
\begin{Def}
	Given a Lie algebra $\mathfrak{g}$, it is called solvable if the derived series: $\mathfrak{g}^{(0)}=\mathfrak{g}, \mathfrak{g}^{(n+1)}:=[\mathfrak{g}^{(n)},\mathfrak{g}^{(n)}]$ terminates; and it is called Nilpotent if the lower central seires: $\mathfrak{g}^0=\mathfrak{g}, \mathfrak{g}^{n+1}:=[\mathfrak{g}^{n},\mathfrak{g}]$ terminates. The unique maximal solvable ideal of $\mathfrak{g}$ is denoted by $Rad(\mathfrak{g})$, and $\mathfrak{g}$ is called semisimple if $Rad(\mathfrak{g})=0$.
\end{Def}
It is obvious that $\mathfrak{g}/Rad(\mathfrak{g})$ is always semisimple.

\begin{Def}
	Given a Lie algebra $\mathfrak{g}$, the Borel subalgebras of $\mathfrak{g}$ are defined to be the maximal solvable subalgebras of $\mathfrak{g}$.
\end{Def}

\begin{Prop} \cite[$\S$4.1, Cor.~4.1,Lie's Theorem]{Humphreys-Lie-Alg} \label{Lie-Thm}
	If $\mathfrak{g}$ is solvable, then with respect to a suitable basis, all elements in $\mathfrak{g}$ are upper triangular. 
\end{Prop}

Recall that every semisimple Lie algebra $\mathfrak{g}$ has its root system \cite[chap.~\RNum{2} to \RNum{4}]{Humphreys-Lie-Alg} such that $\mathfrak{g}=H\oplus \sum\limits_{\alpha \in \Phi} \mathfrak{g}_{\alpha}$, where $H$ is a maximal toral subalgebra of $\mathfrak{g}$ (i.e. a subalgebra consisting of semisimple elements) and $\Phi$ is the set of roots, i.e., nonzero elements in the dual space $H^*$ such that there exists $x\in \mathfrak{g}$, where $ [h,x]=\alpha(h)x$ for all $h\in H$. Moreover, for each $\alpha\in \Phi$, we have a triple $(h_{\alpha}, x_{\alpha}, y_{\alpha})$, where $x_{\alpha}\in {\mathfrak{g}}_{\alpha}, y_{\alpha}\in {\mathfrak{g}}_{-\alpha}, h_{\alpha}=[x_{\alpha},y_{\alpha}]\in H$. 

\begin{Prop}\label{ss-alg-dim5}
	If $\mathfrak{g}$ is semisimple Lie algebra of dimension $\leq 5$, then $\mathfrak{g}$ is simple, and isomorphic to $\mathfrak{sl}_2$ as an abstract Lie algebra.
\end{Prop}
\begin{pf}
	In fact, if $\dim \mathfrak{g}\leq 5$, then the $\Phi$ has at most two elements since otherwise there are at least four roots, which implies $\dim H\geq 2$ and hence $\dim\mathfrak{g}\geq 6$, contradiction. On the other hand, $\dim \sum\limits_{\alpha \in \Phi} \mathfrak{g}_{\alpha}\geq 2$ (since otherwise $\mathfrak{g}=H$. But $H$ is abelian \cite[Lemma.~8.1]{Humphreys-Lie-Alg}, thus solvable, contradiction). Now the proposition follows from the fact that $\mathfrak{sl}_2$ is the only dimensional $3$ semisimple Lie algebra up to isomorphism.
\end{pf}

\bibliographystyle{amsalpha}
\bibliography{mybib}

\end{document}